\documentclass[11pt,a4paper]{article}

\usepackage[T1]{fontenc}
\usepackage[utf8]{inputenc}
\usepackage[american]{babel}
\usepackage{amsmath,amsthm}
\usepackage{txfonts}
\usepackage[margin=2.5cm]{geometry}
\usepackage{microtype}
\usepackage{booktabs}
\usepackage{graphicx}
\AtBeginDocument{\setlength{\unitlength}{0.0025\textwidth}}% picture units scale with the text width
\usepackage{cite}
\usepackage{xcolor}
\usepackage[hidelinks]{hyperref}

\newcommand{\bbR}{\mathbb{R}}
\newcommand{\FF}{\mathcal{F}}
\newcommand{\pa}{\partial}
\newcommand{\dd}{\textrm{d}}
\newcommand{\bbr}{\mathbf{r}}
\newcommand{\bbl}{\boldsymbol{\ell}}
\newcommand{\Jac}{\mathcal{J}}
\newcommand{\cC}{\mathcal{C}}
\newcommand{\muw}{\mu_{\mathrm{w}}}
\newcommand{\mug}{\mu_{\mathrm{g}}}
\newcommand{\muo}{\mu_{\mathrm{o}}}
\newcommand{\muwg}{\mu_{\mathrm{wg}}}
\newcommand{\muwgo}{\mu_{\mathrm{wgo}}}
\newcommand{\wdg}{\wedge}

\theoremstyle{plain}
\newtheorem{theorem}{Theorem}
\newtheorem{proposition}[theorem]{Proposition}
\newtheorem{corollary}[theorem]{Corollary}
\newtheorem{lemma}[theorem]{Lemma}
\theoremstyle{remark}
\newtheorem{remark}[theorem]{Remark}

\newcommand{\preprintlegend}{Preprint submitted to the Journal of Hyperbolic Differential Equations on November 2, 2026.}

\title{On the phase-space geometry for Riemann solutions:\\
rank-one structure, implicit solutions\\
and elusive transitional shocks\thanks{\preprintlegend}}
\author{Pablo Casta\~neda\\[2pt]
\small Departamento Acad\'emico de Matem\'aticas, ITAM, 01080, Cd.\ de M\'exico, Mexico\\
\small Instituto de Matem\'aticas, Ciudad Universitaria, UNAM, 04510, Cd.\ de M\'exico, Mexico\\
\small \texttt{pablo.castaneda@itam.mx}}
\date{}

\begin{document}
\maketitle

\begin{abstract}
 The effective flux function reduces a Riemann problem for a $2\times2$ system
of conservation laws to a scalar one, by lifting the wave curve.
We show that the reduction rests on a rank-one property, namely that the Jacobian
$\pa(u,\,v)/\pa(x,\,t)$ equals the product of the two wave strengths, up to
nonzero factors, and vanishes identically for every self-similar solution.
Hence the lifting needs no hypothesis for Riemann data, and for Cauchy data it
survives exactly on an integral curve of an eigenvector field, where we give an
implicit solution; off that curve the state path leaves it at a rate
we compute. We then ask when an implicit closed-form solution exists for arbitrary
data. For homogeneous quadratic fluxes these are precisely
the systems linearly equivalent to two decoupled Burgers equations (the curve
$a=1+b^2$ of the Schaeffer--Shearer normal form). Every such system loses strict
hyperbolicity along a full line, so no model with an isolated umbilic point
qualifies, among them the quadratic Corey model of three-phase flow. We close
with a quadratic model in which a datum on a straight integral curve gives an
exact solution with a stationary transitional shock, having one family nearly
linearly degenerate at each end.
\end{abstract}



%==============================================================================
\section{Introduction}
\label{sec:intro}

We consider the $2\times2$ system of conservation laws in one spatial
dimension $x\in\bbR$, $t\ge0$,
\begin{equation}
\left\{
\begin{array}{rcl}
 u_t &+& f(U)_x \;=\; 0,\\
 v_t &+& g(U)_x \;=\; 0,
\end{array}
\right.
\qquad U=(u,\,v)^\top\in D\subset\bbR^2,
\label{eq:system}
\end{equation}
with real-valued fluxes $f,\,g$ defined on the phase space $D$, together with
either Riemann data
\begin{equation}
U(x,\,0)=\begin{cases} L=(u_L,\,v_L)^\top, & x<0,\\
                       R=(u_R,\,v_R)^\top, & x>0,\end{cases}
\label{eq:RP}
\end{equation}
denoted RP$(L,\,R)$, or general Cauchy data
\begin{equation}
U(x,\,0)=U_0(x)=(u_0(x),\,v_0(x))^\top .
\label{eq:Cauchy}
\end{equation}
Write $J=\pa F/\pa U$ for $F=(f,\,g)^\top$, and assume strict hyperbolicity
where needed, that is, $J$ has real eigenvalues $\lambda_1<\lambda_2$ with right
eigenvectors $\bbr_1,\,\bbr_2$ and left eigenvectors $\bbl_1,\,\bbl_2$
normalized by $\bbl_i\bbr_j=\delta_{ij}$.

For \eqref{eq:RP} the solution is classically built from a slow wave curve
issuing forward from $L$, in the direction of increasing speeds of the smaller
eigenvalue, and a fast wave curve issuing backward from $R$, in the direction
of decreasing speeds of the larger eigenvalue; the two meet at
an intermediate state $M$, see \cite{CetA16,Lax57,Liu74,Liu75,Rod23}, thus
defining the Riemann solution by the wave curve method. In
\cite{Cas18} the complete wave curve is lifted to a scalar flux, the effective
flux function (EFF), and the solution of the system is recovered from the
scalar Ole\u{\i}nik construction of \cite{Cas16}. The Banana theorem of
\cite{AetA10}, proved there to solve a more general problem, contains the first
implicit use of the EFF, employed to find the separatrix displayed in
\cite{CetA16}. Transitional waves
\cite{Loz18,LetA24} may break the two-group composition; this does not affect
the results below.

Let us state the three questions and the three answers before going further.
The questions are what geometric hypothesis makes the EFF work, how far the
construction reaches beyond Riemann data, and when an implicit closed-form
solution can exist at all. The answers are, in order, (1) the EFF is a
rank-one phenomenon, which for self-similar solutions holds automatically
(Proposition~\ref{prop:rank1} and identity~\eqref{eq:JacAB}); (2) for Cauchy
data the construction is exact on data supported on an integral curve of an eigenvector field
(Theorem~\ref{thm:exact}), and off that set it degrades at a rate controlled
by a single invariant quantity (Proposition~\ref{prop:error}); and (3) a
closed-form implicit solution for arbitrary data exists only for a very thin
class of systems, which for quadratic fluxes we classify completely
(Theorem~\ref{thm:classif}) and which excludes every model with an isolated umbilic point
(Corollary~\ref{cor:nogo}).

In Sec.~\ref{sec:EFF} we recall the lifting and then identify what makes it
possible, namely that a self-similar solution sweeps a curve rather than an
open set in phase space; this is the hypothesis that we must impose by hand in
Sec.~\ref{sec:implicit} once the data are no longer self-similar, which is why
the implicit solution there lives on an integral curve of an eigenvector field and not on an
arbitrary one. In Sec.~\ref{sec:classif} we ask what it would take for an
arbitrary datum to work, and the answer is restrictive enough to be a
classification. In Sec.~\ref{sec:example} we run the whole argument on the
quadratic model, where the numbers are explicit, a shock that crosses the
characteristic speeds (a crossing discontinuity, defined in
Sec.~\ref{sec:transitional}) can be read through invariant quantities, and an
exact solution with a transitional shock can be written down.

%==============================================================================
\section{Lifting, the effective flux function, and rank-one structure}
\label{sec:EFF}

This section has two halves. The first recalls the EFF of \cite{Cas18}, which
turns a Riemann problem for \eqref{eq:system} into the scalar problem
\eqref{eq:scalar} by restricting the fluxes to the wave curve. The second
half asks why a wave curve is available at all, and answers with two
identities, the first being that the Jacobian of the map $(x,\,t)\mapsto U$
vanishes for every self-similar solution, and the second that the same
Jacobian is the product of the strengths of the two waves. The second identity
is used in every later section.

\subsection{The lifting}
\label{sec:lifting}

Assume the solution of \eqref{eq:system}--\eqref{eq:RP} admits a wave curve in
$D$ parametrized as $(u,\,\gamma(u))^\top$, with
$(u_\alpha,\,\gamma(u_\alpha))^\top=\alpha$ for $\alpha=L,\,M,\,R$, where $M$ stands for the intermediate state of the
Riemann solution described in Sec.~\ref{sec:intro}. Then the scalar Riemann
problem
\begin{equation}
u_t+\FF(u)_x=0,\qquad
u(x,\,0)=\begin{cases}u_L,&x<0\\ u_R,&x>0\end{cases}
\qquad\text{with}\qquad \FF(u)=f(u,\,\gamma(u)),
\label{eq:scalar}
\end{equation}
solved by Ole\u{\i}nik's convex hull construction \cite{Cas16}, yields a
solution of \eqref{eq:system}--\eqref{eq:RP} through
$v(x,\,t)=\gamma(u(x,\,t))$. This is Theorem~1 of \cite{Cas18}, and $\FF$ is
the EFF.

Across a shock of speed $\sigma$ joining a left state
$L=(u_L,\,v_L)^\top$ and a right state $R=(u_R,\,v_R)^\top$ of the wave curve,
as in \eqref{eq:RP}, the first Rankine--Hugoniot equation of \eqref{eq:system}
reads $\sigma(u_R-u_L)=f(R)-f(L)$, and since $v=\gamma(u)$ this is
$\sigma(u_R-u_L)=\FF(u_R)-\FF(u_L)$, the Rankine--Hugoniot condition of
\eqref{eq:scalar}. Hence \eqref{eq:scalar} is satisfied along shocks
trivially, and the second equation of \eqref{eq:system} plays no role in this
step.

The quasilinear form $u_t+\FF'(u)u_x=0$ shows that the propagation speed is
$\FF'(u)$, which along a rarefaction coincides with the eigenvalue of the
corresponding family, since:
\begin{equation}
\FF'(u)=\frac{\dd}{\dd u}f(u,\,\gamma(u))=f_u+\gamma'(u)f_v=\lambda .
\label{eq:Fprime}
\end{equation}
Along an integral curve of an eigenvector field of $J$ (a rarefaction curve once it is oriented by increasing eigenvalue; the orientation plays no role below) the eigenvalue problem $J\bbr=\lambda\bbr$ holds
with $\bbr=(1,\,\gamma')^\top$, that is
\begin{equation}
f_u+\gamma'f_v=\lambda,\qquad g_u+\gamma'g_v=\lambda\gamma' .
\label{eq:rarefaction}
\end{equation}
Eliminating $\lambda$ between the two gives the ODE of these integral curves,
\begin{equation}
f_v\,(\gamma')^2+(f_u-g_v)\,\gamma'-g_u=0,
\qquad
\gamma'(u)=\frac{g_v-f_u\pm\sqrt{(g_v-f_u)^2+4f_vg_u}}{2f_v},
\label{eq:gammaPrime}
\end{equation}
and substituting back into the first equation of \eqref{eq:rarefaction}
returns $\lambda=\frac12\big(\operatorname{tr}J\pm
\sqrt{(\operatorname{tr}J)^2-4|J|}\big)$, the eigenvalues of $J$.

Equation \eqref{eq:gammaPrime} breaks down where $f_v(U)=0$, and it is worth
saying what happens there. By \eqref{eq:rarefaction}, $f_v=0$ forces
$\lambda=f_u$ and $g_u=(\lambda-g_v)\gamma'$, so one of the two integral
curves through such a state has $\gamma'=\infty$, so that it is tangent to the
$v$-axis direction and $u$ has a local extremum along it. The state is a
\emph{turning point} of the parametrization by $u$ (not a singularity of the
curve, which remains smooth), and one continues past it by taking $v$ as
parameter, with $\dd u/\dd v=0$ there. Only one of the two branches of
\eqref{eq:gammaPrime} degenerates, the other tending to $g_u/(f_u-g_v)$.

\begin{remark}[quadratic Corey model]
\label{rem:oil}
In petroleum engineering the simplest fluxes form the quadratic Corey model.
Here $u$ and $v$ are the water and gas saturations, $1-u-v$ is the oil
saturation, $\muw,\,\mug,\,\muo$ are the phase viscosities, and $f$ is the
fractional flow of water,
\begin{gather*}
f(U)=\frac{u^2/\muw}{u^2/\muw+v^2/\mug+(1-u-v)^2/\muo},\\[4pt]
\frac{\pa f}{\pa v}
=\frac{2\,(u^2/\muw)\big((1-u-v)/\muo-v/\mug\big)}
      {\big(u^2/\muw+v^2/\mug+(1-u-v)^2/\muo\big)^2},
\end{gather*}
so that $f_v$ vanishes only for $u=0$ or $v\muo=(1-u-v)\mug$. These are two
straight lines in the saturation triangle, characterized for this purpose in
\cite{AetA10}. By the previous paragraph they are exactly where the
parametrization of the wave curve by water saturation must be exchanged for
one by gas saturation, so the turning points of the EFF construction are known
in closed form for this model, and the same identification carries over to the
more intricate models of \cite{CMF22}.

Those lines are not the only ones that matter. Let $W=(1,\,0)$, $G=(0,\,1)$ and $O=(0,\,0)$ be the
vertices of the triangle where the water, gas and oil saturations equal one,
and let $B$ be the point of the side $WG$ with $s_g=(\mug/\muw)\,s_w$. The
line $BO$, $u/\muw=v/\mug$, contains the umbilic point
$U=(\muw,\,\mug)/(\muw+\mug+\muo)$, where $u/\muw=v/\mug=(1-u-v)/\muo$, and
along it $f/g=\muw/\mug$ is constant, so that $F(U)$ is parallel to
$(\muw,\,\mug)$ there. Hence $F$ maps the line into itself, any two of its
states satisfy the Rankine--Hugoniot condition, and $(\muw,\,\mug)$ is an
eigenvector of $J$ along the line, which is therefore an integral curve of an
eigenvector field (Remark~\ref{rem:corey}). We follow the notation of
\cite{AetA10}.
\end{remark}

\subsection{Rank one, or why the lifting exists}
\label{sec:rankone}

That the state path of the solution is a curve $v=\gamma(u)$ is a hypothesis, and for Cauchy data it can fail for an unremarkable physical
reason, namely that along a profile the water saturation may well increase and
then decrease, so $u_0$ is not invertible and $v_0\circ u_0^{-1}$ defines no
function. For Riemann data, however, the hypothesis costs nothing, at least on
each segment of the wave curve between turning points, for the reason given in
the paragraph after \eqref{eq:gammaPrime}.

\begin{proposition}
\label{prop:rank1}
Let $U(x,\,t)=W(x/t)$ be any self-similar solution of \eqref{eq:system} with
$W=(W_1,\,W_2)^\top$ differentiable. Then
\[
\Jac:=\frac{\pa(u,\,v)}{\pa(x,\,t)}=u_xv_t-u_tv_x\equiv 0 .
\]
\end{proposition}

\begin{proof}
With $\xi=x/t$ one has $u_x=W_1'/t$, $u_t=-\xi W_1'/t$, $v_x=W_2'/t$ and
$v_t=-\xi W_2'/t$, whence $u_xv_t-u_tv_x=(-\xi W_1'W_2'+\xi
W_1'W_2')/t^2$ vanishes, as claimed.
\end{proof}

Thus for Riemann data the map $(x,\,t)\mapsto U(x,\,t)$ has rank one and its
image is a curve, the wave curve. The EFF is the lifting of that curve, and
the whole of \cite{Cas18} is a statement about rank-one solutions. Let us now
record the identity that explains the mechanism, and that will be used
repeatedly.

\begin{proposition}
\label{prop:JacAB}
For any $\cC^1$ solution of \eqref{eq:system} in a strictly hyperbolic region,
set the wave strengths
\[
\kappa_1:=\bbl_1\cdot U_x,\qquad \kappa_2:=\bbl_2\cdot U_x .
\]
Then
\begin{equation}
\Jac \;=\; -\,(\lambda_2-\lambda_1)\,\det[\bbr_1,\,\bbr_2]\;\kappa_1\kappa_2 .
\label{eq:JacAB}
\end{equation}
\end{proposition}

\begin{proof}
Write $U_x=\kappa_1\bbr_1+\kappa_2\bbr_2$, which is legitimate because
$\{\bbr_1,\,\bbr_2\}$ is a basis and $\{\bbl_1,\,\bbl_2\}$ its dual (from
$\lambda_1(\bbl_1\bbr_2)=(\bbl_1J)\bbr_2=\bbl_1(J\bbr_2)=
\lambda_2(\bbl_1\bbr_2)$ and $\lambda_1\neq\lambda_2$ we get
$\bbl_1\bbr_2=0$). From \eqref{eq:system},
$U_t=-JU_x=-(\lambda_1\kappa_1\bbr_1+\lambda_2\kappa_2\bbr_2)$. Hence
\[
\Jac=\det[U_x,\,U_t]
=-\det\big[\kappa_1\bbr_1+\kappa_2\bbr_2,\;\lambda_1\kappa_1\bbr_1+\lambda_2\kappa_2\bbr_2\big]
=-(\lambda_2-\lambda_1)\,\kappa_1\kappa_2\,\det[\bbr_1,\,\bbr_2],
\]
which completes the proof.
\end{proof}

The content of \eqref{eq:JacAB} is that the equation forces $U_t$ to be the
image of $U_x$ under $J$, so that $\det[U_x,\,U_t]$ measures how far $U_x$ is
from being an eigenvector. If $\kappa_2=0$ then $U_t=-\lambda_1U_x$ is parallel to
$U_x$, and the Jacobian vanishes by collinearity rather than by cancellation
between terms. The identity does not depend on how the eigenvectors are
scaled, since $\bbr_i\mapsto c_i\bbr_i$ sends $\kappa_1\mapsto\kappa_1/c_1$,
$\kappa_2\mapsto\kappa_2/c_2$ and $\det[\bbr_1,\,\bbr_2]\mapsto c_1c_2\det[\bbr_1,\,\bbr_2]$.

\begin{remark}
\label{rem:threefold}
Three consequences will be used without further comment. First, $\Jac=0$ if
and only if $\kappa_1\kappa_2=0$, that is, if and only if one of the two families carries no
strength; this is the definition of a simple wave, so rank one, the existence
of $\gamma$, and the coincidence of the two component speeds of
Sec.~\ref{sec:implicit} are one and the same condition. Second, strict
hyperbolicity enters \eqref{eq:JacAB} through the factor
$\lambda_2-\lambda_1$, rather than as a separate hypothesis. Third, combining
with Proposition~\ref{prop:rank1}, in a Riemann solution $\kappa_1\kappa_2\equiv0$ holds, so at
almost every point at most one family is active. That is why the state path of
a Riemann solution is a curve, and it is the structural reason behind the EFF.
\end{remark}

%==============================================================================
\section{Implicit solutions for Cauchy data}
\label{sec:implicit}

A scalar conservation law may be seen as a crowd walking in single file, in
which every state carries its own speed and nothing meets anything else until a
shock forms.
For the single file the implicit solution is classical, and although it
defines the solution implicitly, it is a closed formula. A system is two such
crowds sharing one corridor, each walking at a speed that depends on both, and
in general there is no such formula. Treating the two crowds as independent
means transporting each component at its own speed, and the two speeds may
differ. This section asks when the crowds may still be treated as
independent, and how fast they stop being so when they may not.
In Sec.~\ref{sec:template} we define the two component speeds and show that they
coincide exactly in the rank-one case of Sec.~\ref{sec:rankone}; in
Sec.~\ref{sec:exactresult} we turn the coincidence into a theorem; in
Sec.~\ref{sec:firstorder} we measure the error when it fails; and in
Sec.~\ref{sec:exactrep} we record what remains true with no hypothesis at all.

\subsection{The scalar template and the two component speeds}
\label{sec:template}

For a scalar problem $u_t+f(u)_x=0$ with $u(x,\,0)=u_0(x)$ smooth, the
implicit solution is $u(x,\,t)=u_0(\xi)$ with $\xi=x-t\lambda$ and
$\lambda=f'(u(x,\,t))$, valid for
$t<\big(\sup_\xi[-f''(u_0(\xi))u_0'(\xi)]\big)^{-1}$, or for all $t$ if the
supremum is nonpositive. What makes the formula \emph{closed} is that $f'(u)$
is constant along each characteristic, so the characteristic is a straight
line and $x-t\lambda$ is its foot at $t=0$.

For a function $U(x,\,t)$ we write $\Gamma_t:=\{U(x,\,t):\,x\in\bbR\}$ for its
\textbf{state path} at time $t$, the curve traced in phase space by the
profile at that instant. Thus $\Gamma_0$ is the curve drawn by the datum, and
$\Gamma_t$ is independent of how the profile is parametrized.

Now let $u_0,\,v_0$ be monotone and differentiable, so that $v_0=\gamma\circ
u_0$ defines $\gamma$, and set
\begin{equation}
\lambda:=f_u+\gamma'(u)f_v,\qquad
\mu:=\frac{1}{\gamma'(u)}g_u+g_v ,
\label{eq:lammu}
\end{equation}
which are exactly the factors making $\pa_xf=\lambda u_x$ and
$\pa_xg=\mu v_x$. Equivalently, for any solution,
\begin{equation}
\lambda=-\frac{u_t}{u_x},\qquad \mu=-\frac{v_t}{v_x},
\qquad\text{and}\qquad
\Jac=u_xv_x\,(\lambda-\mu).
\label{eq:lammu2}
\end{equation}
So $\lambda$ and $\mu$ are the slopes $\dd x/\dd t$ of the level curves of $u$
and of $v$, and by \eqref{eq:JacAB} and \eqref{eq:lammu2},
\begin{equation}
\lambda=\mu
\iff \Jac=0
\iff \kappa_1\kappa_2=0
\iff (1,\,\gamma')^\top \text{ is an eigenvector of } J .
\label{eq:equiv}
\end{equation}

\begin{remark}
\label{rem:noninvariance}
The quantities $\lambda$ and $\mu$ are intrinsic to the functions $u$ and $v$
but not to the system, since they depend on which pair of functions is used as state
variables. In the model of Sec.~\ref{sec:example}, at $U=(1,\,1)$ with
$u_x=v_x$ one gets $\lambda=7$ and $\mu=4$. In the variables
$(u,\,u-\frac12v)$, which are equally monotone and still conservative, one
gets $\lambda=7$ and $\mu=10$, reversing the order; in $(u,\,u+v)$ one gets
$\mu=\frac{11}{2}$, and in $(u,\,u+2v)$ one gets $\mu=5$. What is invariant is the equality itself, since by \eqref{eq:equiv} $\lambda=\mu$
holds if and only if $(1,\,\gamma')^\top$ is an eigenvector of $J$, and every
change of dependent variables sends eigenvectors to eigenvectors; only the
order of $\lambda$ and $\mu$, and the size of their difference, depend on the
variables. Admissibility
criteria are invariant under changes of dependent variables, so no criterion
phrased in terms of $\lambda$ and $\mu$ can be equivalent to one. This is why
the discussion in Sec.~\ref{sec:example} is carried out with
$\nabla\lambda_i\cdot\bbr_i$ instead.
\end{remark}

\subsection{The exact result}
\label{sec:exactresult}

Why look for an implicit solution of a system at all, when the Riemann problem
is already solved by wave curves? Because the Riemann solution is a single
self-similar object, whereas a reservoir, or a road, starts from a smooth
profile and the question is how that profile deforms before any wave forms.
For a scalar law the implicit formula answers this at once, and it is what
underlies the Lagrangian--Eulerian schemes of \cite{AetA22,AetA23}.
The following theorem says how much of that survives for a system, namely all of it,
provided the initial profile is drawn on an integral curve of an eigenvector field.

\begin{theorem}
\label{thm:exact}
Let $f,\,g\in \cC^2(D)$ and let $\gamma$ solve \eqref{eq:gammaPrime}, that is,
let $\Gamma=\{(u,\,\gamma(u))\}$ be an integral curve of one of the
eigenvector fields of $J$, with $f_v\neq0$ along it. Let $u_0\in \cC^1(\bbR)$
with $(u_0,\,\gamma\circ u_0)\in D$, and take Cauchy data
$U_0=(u_0,\,\gamma\circ u_0)^\top$. Put $\FF(u)=f(u,\,\gamma(u))$ and
$\lambda(u)=\FF'(u)=f_u+\gamma'f_v$. Then for
\[
0\le t<t_b:=\Big(\sup_\xi\big[-\lambda'(u_0(\xi))\,u_0'(\xi)\big]\Big)^{-1}
\]
the relation $u=u_0\big(x-t\lambda(u)\big)$ defines a unique $\cC^1$ function,
and
\[
U(x,\,t)=\big(u(x,\,t),\,\gamma(u(x,\,t))\big)
\]
solves \eqref{eq:system}. Moreover $\mu\equiv\lambda$, so the companion
relation $v=v_0(x-t\mu(v))$ holds automatically and the state
path of Sec.~\ref{sec:template} satisfies $\Gamma_t\equiv\Gamma$ for all
$t<t_b$.
\end{theorem}

\begin{proof}
Since $\Gamma$ is an integral curve of an eigenvector field,
\eqref{eq:rarefaction} holds along it, that is, $f_u+\gamma'f_v=\lambda$ and
$g_u+\gamma'g_v=\lambda\gamma'$. The second gives
$\mu=g_u/\gamma'+g_v=\lambda$. By the scalar theory,
$u=u_0(x-t\lambda(u))$ defines a $\cC^1$ function for $t<t_b$, with
\[
u_x=\frac{u_0'(\xi)}{1+t\,u_0'(\xi)\lambda'(u)},\qquad
u_t=\frac{-\lambda(u)\,u_0'(\xi)}{1+t\,u_0'(\xi)\lambda'(u)},
\]
so that $u_t+\lambda(u)u_x=0$. Setting $v=\gamma(u)$ gives $v_x=\gamma'u_x$
and $v_t=\gamma'u_t$, hence
\[
u_t+f(U)_x=u_t+(f_u+\gamma'f_v)u_x=u_t+\lambda u_x=0,
\]
\[
v_t+g(U)_x=\gamma'u_t+(g_u+\gamma'g_v)u_x
=\gamma'u_t+\lambda\gamma'u_x=\gamma'(u_t+\lambda u_x)=0,
\]
so both equations of \eqref{eq:system} hold, which completes the proof.
\end{proof}

Theorem~\ref{thm:exact} is the Cauchy counterpart of the EFF construction. By
Proposition~\ref{prop:rank1} the Riemann case is rank one for free, and here
rank one is bought with the hypothesis on $\gamma$; the two statements are the
same theorem in the two regimes. In physical terms the hypothesis says that
the initial profile must already be a simple wave, that is, the proportions of the
phases along the profile are those of one characteristic family, and then each
state travels at its own speed along a straight characteristic, exactly as in
the scalar case. The breaking time $t_b$ is the scalar one for the EFF.

\subsection{A first-order estimate off the wave curve}
\label{sec:firstorder}

If the data are not supported on an integral curve of an eigenvector field, the pair
\begin{equation}
u=u_0\big(x-t\lambda(u)\big),\qquad v=v_0\big(x-t\mu(v)\big)
\label{eq:ansatz}
\end{equation}
is not a solution of \eqref{eq:system}, since each relation solves its own
scalar transport equation but the two together solve the system only if the state remains on $\Gamma$ for all times, which by \eqref{eq:equiv} requires $\lambda=\mu$.
What survives is the following.

\begin{proposition}
\label{prop:error}
Let $f,\,g\in\cC^3(D)$ and $u_0,\,v_0\in\cC^3(\bbR)$ with $u_0'$ and $v_0'$
nowhere zero, so that $\gamma=v_0\circ u_0^{-1}$ is defined and $\lambda,\,\mu$
are given by \eqref{eq:lammu}. Let $(u,\,v)$ be defined by \eqref{eq:ansatz}
for
\[
0\le t<\min\Big\{\big(\sup_\xi[-\lambda'(u_0)u_0']\big)^{-1},\
                 \big(\sup_\zeta[-\mu'(v_0)v_0']\big)^{-1}\Big\},
\]
and set $\delta:=v-\gamma(u)$. Then, uniformly for $x$ in compact sets,
\[
\delta(x,\,t)=t\,E(x)+O(t^2),\qquad
E:=u_0'\,\gamma'(u_0)\,\big[\lambda(u_0)-\mu(v_0)\big],
\]
and the residuals $R_1:=u_t+f(U)_x$ and $R_2:=v_t+g(U)_x$ of
\eqref{eq:system} are $O(t)$. In particular \eqref{eq:ansatz} solves
\eqref{eq:system} if and only if $\lambda\equiv\mu$, that is if and only if
$(1,\,\gamma')^\top$ is an eigenvector of $J$.
\end{proposition}

\begin{proof}
Under these hypotheses $\gamma\in\cC^3$, so $\lambda,\,\mu\in\cC^2$, and the
implicit function theorem applied to \eqref{eq:ansatz} gives
$u,\,v,\,\delta\in\cC^2$ in $(x,\,t)$. From Theorem~\ref{thm:exact}, each
relation in \eqref{eq:ansatz} satisfies its own scalar law, so
$u_t=-\lambda(u)u_x$ and $v_t=-\mu(v)v_x$ in the stated range. Hence
$\delta_t=v_t-\gamma'(u)u_t=-\mu v_x+\gamma'\lambda u_x$, and at $t=0$ one has
$v_x=\gamma'u_x$ with $u_x=u_0'$, so that $\delta_t|_{t=0}=E$; since
$\delta(x,\,0)=0$ and $\delta\in\cC^2$, Taylor's formula gives the expansion.
For the residuals, $R_1$ and $R_2$ involve first derivatives of $u$ and $v$,
hence are of class $\cC^1$, and they vanish at $t=0$, since there
$f(U)_x=(f_u+\gamma'f_v)u_0'=\lambda u_0'=-u_t$, and likewise
$g(U)_x=(g_u/\gamma'+g_v)v_0'=\mu v_0'=-v_t$. A $\cC^1$ function that vanishes
at $t=0$ is $O(t)$ on compact sets. The last assertion is \eqref{eq:equiv}.
\end{proof}

It is curious that, to first order, the implicit construction works for every
smooth datum, with no hypothesis on $\gamma$. The residuals vanish at $t=0$, so
the predictor has the same initial velocity $U_t(x,\,0)=-F(U_0)_x$ as the
classical solution (which exists for short times when the system is strictly
hyperbolic along the data), whether or not $\lambda=\mu$, and the two agree to
first order in $t$. The hypothesis $\lambda\equiv\mu$ of
Theorem~\ref{thm:exact} is what makes the agreement exact. (The term $tE$ in
$\delta$ is then not an error of the predictor, since the classical solution leaves
$\Gamma$ with the same $E$.)

\begin{remark}
\label{rem:lowreg}
The hypotheses of Proposition~\ref{prop:error} can be relaxed at the price of the rates. With
$f,\,g\in\cC^2(D)$ and $u_0,\,v_0\in\cC^2(\bbR)$, which is what
\eqref{eq:lammu} needs in order to make sense, the same argument gives
$\delta\in\cC^1$ and $R_1,\,R_2\in\cC^0$, hence $\delta=t\,E+o(t)$ and
$R_1,\,R_2=o(1)$ as $t\to0$, still uniformly on compact sets. The leading term
$E$ survives; only the size of the remainder is lost.
\end{remark}

Thus \eqref{eq:ansatz} is a first-order predictor whose error is governed by
the invariant obstruction $|\lambda-\mu|$, equivalently by $|\kappa_1\kappa_2|$ through
\eqref{eq:JacAB} and \eqref{eq:lammu2}. This is the precise sense in which the
two component speeds ``remain close for small times'', and it is the form in
which the construction is useful for the schemes of
\cite{AetA22,AetA23}, since the initial speeds are explicit from the data
and Proposition~\ref{prop:error} says how long one may propagate each phase at
its own speed before the two profiles have to be reconciled.

\subsection{The deformation of the state path}
\label{sec:deformation}

\begin{figure}[htbp]
\centering
\small
\begin{picture}(200,145)(0,-16)
\put(0, 0){\includegraphics[width=0.5\linewidth]{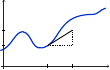}}

\put(190, -5){$u$}
\put( 89, -5){\makebox(0,0)[tc]{$u(x,\,t)$}}
\put(134, -16){\makebox(0,0)[tc]{$u(x,\,t) + \Delta u$}}

\put(-1, 122){$v$}
\put(  3, 36){\makebox(0,0)[rb]{$\gamma(u)$}}
\put(  3, 64){\makebox(0,0)[rb]{$\gamma(u) + \Delta v$}}

\put(97, 33){$1 \approx \Delta u$}
\put(140, 53){$\gamma'(u) \approx \Delta v$}

\put(160, 108){\color{blue} $\gamma(u)$}
\end{picture}
\caption{In the limit, the slope of the state path is the ratio
$v_x(x,\,t)/u_x(x,\,t)$ of \eqref{eq:gammat}. Its explicit form is
\eqref{eq:gammatExplicit}, where the two spatial derivatives are evaluated at the
different arguments $\xi=x-t\lambda(u)$ and $\zeta=x-t\mu(v)$, and each
carries its own denominator. Hence the path deforms in time unless
$\lambda\equiv\mu$.}
\label{fig:gamma}
\end{figure}

Let us record how the state path itself moves. At time $t$ the state
path $\Gamma_t$ of the pair \eqref{eq:ansatz}, the image of
$x\mapsto(u(x,\,t),\,v(x,\,t))$, has slope at the point of parameter $x$ equal
to the ratio of the spatial derivatives,
\begin{equation}
\gamma_t'(u)=\frac{v_x(x,\,t)}{u_x(x,\,t)} ;
\label{eq:gammat}
\end{equation}
see Fig.~\ref{fig:gamma}. Since $u=u_0(\xi)$ with $\xi=x-t\lambda(u)$ and
$v=v_0(\zeta)$ with $\zeta=x-t\mu(v)$, differentiation gives
$u_x=u_0'(\xi)/\big(1+t\,u_0'(\xi)\lambda'\big)$ and
$v_x=v_0'(\zeta)/\big(1+t\,v_0'(\zeta)\mu'\big)$, so that
\begin{equation}
\gamma_t'(u)=\frac{v_0'(\zeta)}{u_0'(\xi)}\cdot
\frac{1+t\,u_0'(\xi)\,\lambda'(u)}{1+t\,v_0'(\zeta)\,\mu'(v)} .
\label{eq:gammatExplicit}
\end{equation}
At $t=0$ this reduces to $\gamma'(u)$, as it must. The two arguments
$\xi\neq\zeta$ and the two denominators are what makes $\Gamma_t$ deform, and
they coincide precisely when $\lambda=\mu$, in agreement with
\eqref{eq:equiv} and Theorem~\ref{thm:exact}.

\subsection{The exact representation for general data}
\label{sec:exactrep}

For arbitrary Cauchy data the correct statement uses Riemann invariants
$w,\,z$, characterized by $\nabla w=\bbl_1$ and $\nabla z=\bbl_2$, which
satisfy $w_t+\lambda_1w_x=0$ and $z_t+\lambda_2z_x=0$. Then, always,
\begin{equation}
w(x,\,t)=w_0\big(X_1(x,\,t)\big),\qquad z(x,\,t)=z_0\big(X_2(x,\,t)\big),
\label{eq:exactrep}
\end{equation}
where $X_i$ is the foot at $t=0$ of the $i$-characteristic, determined by
$\pa_tX_i+\lambda_i\pa_xX_i=0$ and $X_i(x,\,0)=x$. This has the right number
of arbitrary functions and requires no hypothesis, but $X_i$ is defined by an
ODE rather than by a formula. The substitution $X_i=x-t\lambda_i$ is
legitimate precisely when $\lambda_i$ is constant along its own
characteristic, that is when
\begin{equation}
\lambda_1=\lambda_1(w)\qquad\text{and}\qquad \lambda_2=\lambda_2(z).
\label{eq:decoupling}
\end{equation}
Sec.~\ref{sec:classif} determines when \eqref{eq:decoupling} holds.

\begin{remark}
\label{rem:gas}
The classical instance of \eqref{eq:decoupling} is isentropic gas dynamics,
\[
\rho_t+(\rho\,\hat u)_x=0,\qquad (\rho\,\hat u)_t+\big(\rho\,\hat u^2
+p(\rho)\big)_x=0,\qquad p(\rho)=K\rho^{\,3},
\]
where $\rho$ is the density, $\hat u$ the velocity and $p$ the pressure, which
here is proportional to the cube of the density. Then the speed of sound
$c=\sqrt{p'(\rho)}$ is proportional to $\rho$, and therefore the Riemann
invariants are $r=\hat u-c$ and $s=\hat u+c$, while the eigenvalues are
$\lambda_1=r$ and $\lambda_2=s$, so that the eigenvalues are the Riemann
invariants and the system decouples into two Burgers equations \cite{CF48}. See
Remark~\ref{rem:gamma3} for the role this example plays in delimiting
Theorem~\ref{thm:classif}.
\end{remark}

%==============================================================================
\section{When does the implicit formula work for every datum?}
\label{sec:classif}

Define, for a strictly hyperbolic system,
\begin{equation}
A:=\nabla\lambda_1\cdot\bbr_2,\qquad B:=\nabla\lambda_2\cdot\bbr_1 .
\label{eq:AB}
\end{equation}
Condition \eqref{eq:decoupling} is exactly $A\equiv B\equiv0$. Let us first
check that this criterion, unlike the component speeds of
Remark~\ref{rem:noninvariance}, does not depend on the variables in which the
model happens to be written.

\begin{lemma}[invariance]
\label{lem:invar}
Let $\Psi$ be a $\cC^1$ diffeomorphism of the dependent variables,
$\widetilde U=\Psi(U)$, and let $\widetilde A$ and $\widetilde B$ be the
quantities \eqref{eq:AB} of the system written in the variables
$\widetilde U$, computed with the eigenvectors
$\widetilde\bbr_i=D\Psi\,\bbr_i$. Then $\widetilde A=A$ and
$\widetilde B=B$ at corresponding points. In particular $A\equiv B\equiv0$
in one set of variables if and only if it holds in the other.
\end{lemma}

\begin{proof}
From $\widetilde U_t=D\Psi\,U_t=-D\Psi\,J\,U_x
=-\big(D\Psi\,J\,(D\Psi)^{-1}\big)\widetilde U_x$ the system in the variables
$\widetilde U$ has matrix $\widetilde J=D\Psi\,J\,(D\Psi)^{-1}$, which is
conjugate to $J$. Hence $\widetilde\lambda_i\circ\Psi=\lambda_i$, the vectors
$\widetilde\bbr_i=D\Psi\,\bbr_i$ are right eigenvectors of $\widetilde J$, and
the chain rule gives
$\widetilde\nabla\widetilde\lambda_i=(D\Psi)^{-\top}\nabla\lambda_i$. Therefore
\[
\widetilde\nabla\widetilde\lambda_1\cdot\widetilde\bbr_2
=\nabla\lambda_1^{\top}(D\Psi)^{-1}D\Psi\,\bbr_2=\nabla\lambda_1\cdot\bbr_2 ,
\]
and the same computation applies to $B$.
\end{proof}

For a nonlinear $\Psi$ the transformed system is no longer in conservation
form, which does not matter here, since $A$ and $B$ are built from $J$ alone.
Two
qualifications are needed. The value of $A$ is invariant only with the induced
normalization $\widetilde\bbr_i=D\Psi\,\bbr_i$, since $A$ is linear in
$\bbr_2$, rescaling the eigenvector rescales $A$, and a non-orthogonal
$D\Psi$ does not send unit eigenvectors to unit eigenvectors. What survives
every choice of normalization is whether $A$ vanishes, and that is all this
section uses. The second qualification is a consequence. The condition $A\equiv B\equiv0$ is
invariant under every diffeomorphism, whereas the constancy of the eigendirections
and the decoupling into Burgers equations are invariant only under linear
changes of variables, so the passage from the former to the latter cannot be
expected to hold without a linear structure. A named
instance of the invariance is relativistic isentropic gas dynamics, for which
$A\equiv B\equiv0$ holds in the variables $(\rho,\,u)$ if and only if it holds
in the variables $(\rho_{\mathrm{re}},\,u_{\mathrm{re}})$ of
\cite[Prop.~1.1]{FP24}. With $p=p(\rho)$ the pressure and $c$ the speed of
light, these are
\[
\rho_{\mathrm{re}}=\frac{pu^2+c^4\rho}{c^2(c^2-u^2)},\qquad
u_{\mathrm{re}}=\frac{c^2u\,(p+\rho c^2)}{pu^2+c^4\rho}
\]
and the proposition asserts that
$(\rho,\,u)\mapsto(\rho_{\mathrm{re}},\,u_{\mathrm{re}})$ is a one-to-one
local diffeomorphism on the physical domain $\rho>0$, $0<p'(\rho)<c^2$,
$|u|<c$. It extends \cite[Prop.~1]{ST93} (see also \cite{Chen95}), where the
result originates.

This is what will allow us to reduce any quadratic model to a normal form
below without losing the property under study. The criterion also has a
geometric reading, which is the one to keep in mind.

\begin{lemma}[geometric meaning]
\label{lem:geom}
$A\equiv0$ if and only if $\lambda_1$ is constant along every integral curve
of $\bbr_2$, that is, if and only if $\lambda_1=\lambda_1(w)$ is a function of the
Riemann invariant $w$ of its own family; likewise $B\equiv0$ if and only if
$\lambda_2=\lambda_2(z)$.
\end{lemma}

\begin{proof}
$A=0$ means $\nabla\lambda_1\perp\bbr_2$, which says that $\lambda_1$ is
constant along the integral curves of $\bbr_2$. In the plane the orthogonal
complement of $\bbr_2$ is spanned by $\bbl_1=\nabla w$, since $\bbl_1\bbr_2=0$, so
$\nabla\lambda_1$ is parallel to $\nabla w$, that is, $\lambda_1$ is locally a
function of $w$. The same argument applies to $B$ and $z$.
\end{proof}

So the test is visual, since the integral curves of one family must be the
level curves of the speed of the other. Let us now settle the quadratic case,
which is the local model of a system near an umbilic point in the sense of
\cite{SS87,IMPT88}, and is obtained by translating the umbilic state to the
origin, subtracting its double eigenvalue from the frame velocity, and keeping
the quadratic part of the Taylor expansion of the flux. Write
\begin{equation}
f=\tfrac12\big(au^2+2buv+cv^2\big),\qquad
g=\tfrac12\big(du^2+2euv+kv^2\big),
\label{eq:quadflux}
\end{equation}
so that $f_u=au+bv$, $f_v=bu+cv$, $g_u=du+ev$, $g_v=eu+kv$ are linear forms
and $J$ has linear entries.

For two linear forms $\varphi=\varphi_1u+\varphi_2v$ and
$\psi=\psi_1u+\psi_2v$ write
\begin{equation}
\varphi\wdg\psi:=\varphi_1\psi_2-\varphi_2\psi_1 ,
\label{eq:wedge}
\end{equation}
the determinant of their coefficient vectors in the basis $\{u,\,v\}$; it
vanishes exactly when $\varphi$ and $\psi$ are proportional.

\begin{theorem}[classification]
\label{thm:classif}
Assume that the system with flux \eqref{eq:quadflux} is strictly hyperbolic on
some open set. Then the following are equivalent.
\begin{enumerate}
\item[\textup{(i)}] $A\equiv B\equiv0$ on an open set where the system is
strictly hyperbolic;
\item[\textup{(ii)}] the linear forms $f_v$, $f_u-g_v$ and $g_u$ are pairwise
proportional;
\item[\textup{(iii)}] the matrix
$\mathcal{M}=\begin{pmatrix} b & c\\ a-e & b-k\\ -d & -e\end{pmatrix}$,
whose rows are the coefficient vectors of $f_v$, $f_u-g_v$ and $-g_u$, has
rank $\le1$;
\item[\textup{(iv)}] the eigendirection fields are constant, \emph{i.e.},
$\bbr_1,\,\bbr_2$ do not depend on $U$;
\item[\textup{(v)}] the system is linearly equivalent to two decoupled Burgers
equations, $w_t+\big(\frac12\alpha w^2\big)_x=0$ and
$z_t+\big(\frac12\beta z^2\big)_x=0$ with $w,\,z$ linear in $(u,\,v)$.
\end{enumerate}
\end{theorem}

\begin{proof}
\emph{Elimination of the radical.} Let $m$ be the eigenvector slope,
$\bbr=(1,\,m)^\top$. By \eqref{eq:gammaPrime}, its roots satisfy
$f_vm^2+(f_u-g_v)m-g_u=0$, whose coefficients are linear in $U$, so that
$m_1+m_2=-(f_u-g_v)/f_v$ and $m_1m_2=-g_u/f_v$ are rational in $U$. Swapping
the labels $1\leftrightarrow2$ swaps $A$ and $B$, hence $S:=A+B$ and
$P:=A\,B$ are symmetric in $(m_1,\,m_2)$ and therefore rational in $U$, free
of radicals. Since $A\equiv B\equiv0$ is equivalent to $S\equiv0$ and
$P\equiv0$, the problem is algebraic. Moreover, with the eigenvectors $\bbr_i=(f_v,\,\lambda_i-f_u)^\top$ used
below, $A$ and $B$ are homogeneous of degree $1$, so vanishing on an open set is vanishing on an arc of the unit
circle, and the numerators of $S$ and $P$ are binary forms with polynomial
coefficients in $(a,\dots,k)$.

\emph{Computation.} With the normalization
$\bbr_i=(f_v,\,\lambda_i-f_u)^\top$ one finds
\begin{equation}
A+B=u\,S_1+v\,S_2,\qquad
S_1=f_v\wdg g_u=be-cd,\qquad
S_2=f_v\wdg(g_v-f_u)=ac-b^2+bk-ce,
\label{eq:S}
\end{equation}
The product $P=A\,B$ is not needed in the proof; it is computed in
Remark~\ref{rem:P}.

\emph{Let us go by parts.}

\emph{$(i)\,\Rightarrow\,(iii)$.} Suppose first $f_v\not\equiv0$. From
$A\equiv B\equiv0$ we get $S\equiv0$, hence $S_1=S_2=0$ by \eqref{eq:S}. The
three $2\times2$ minors of $\mathcal M$ are $-S_2$, $-S_1$ and $M_3$; when the
first row $(b,\,c)$ of $\mathcal M$ is nonzero, the vanishing of the first two
forces the third, so that rank $\mathcal M\le1$. If instead $f_v\equiv0$, that
is $b=c=0$, then $J$ is triangular and a direct computation gives the single
condition $e^2-ae-dk=0$, which is $M_3=0$; again rank $\mathcal M\le1$. The
converse is the chain \emph{$(iii)\,\Rightarrow\,(ii)\,\Rightarrow\,(iv)\,\Rightarrow\,(v)\,\Rightarrow\,(i)$}
below.

\emph{$(iii)\,\Longleftrightarrow\,(ii)$.} The rows of $\mathcal M$ are the coefficient vectors
of $f_v$, $f_u-g_v$ and $-g_u$, and rank $\le1$ is pairwise proportionality.

\emph{$(ii)\,\Rightarrow\,(iv)$.} If the three forms are multiples of a common
linear form $\ell$, say $f_v=p\ell$, $f_u-g_v=-q\ell$ and $g_u=r\ell$, the
eigenvector equation factors as $\ell\cdot(pm^2-qm-r)=0$, so where
$\ell\neq0$ the slopes $m_{1,2}$ are the roots of a polynomial with constant
coefficients. These roots are real and distinct, because the system is
strictly hyperbolic on some open set and a polynomial with constant
coefficients has the same roots everywhere. (Without the standing hypothesis
$(ii)$ could hold with complex roots, as for $p\,m^2+p=0$, and the system would
be elliptic; this is why the hypothesis is stated.)

\emph{$(iv)\,\Rightarrow\,(v)$.} Pass to linear coordinates $(w,\,z)$ aligned with
the constant eigendirections. There $J$ is diagonal for all $U$, that is
$\pa F_1/\pa z\equiv0$ and $\pa F_2/\pa w\equiv0$, so $F_1=F_1(w)$ and
$F_2=F_2(z)$; being homogeneous quadratic, $F_1=\frac12\alpha w^2$ and
$F_2=\frac12\beta z^2$.

\emph{$(v)\,\Rightarrow\,(i)$.} In the decoupled pair $\lambda_1=\alpha w$ depends
only on $w$ and $\lambda_2=\beta z$ only on $z$, while $\bbr_1$ points along the
$w$-axis and $\bbr_2$ along the $z$-axis. Hence $A$ is a multiple of
$\pa\lambda_1/\pa z=0$ and $B$ a multiple of $\pa\lambda_2/\pa w=0$, so $A=B=0$,
and by Lemma~\ref{lem:invar} this persists in any linear coordinates. The chain
\emph{$(iii)\,\Rightarrow\,(ii)\,\Rightarrow\,(iv)\,\Rightarrow\,(v)\,\Rightarrow\,(i)$},
together with $(i)\,\Rightarrow\,(iii)$, shows that the five conditions are
equivalent, thus proving the theorem.
\end{proof}

\begin{remark}[the product $AB$]
\label{rem:P}
The proof uses $A+B$ only, but the product has a form worth recording, because
it explains why the Hugoniot locus of the origin enters. With the normalization
$\bbr_i=(f_v,\,\lambda_i-f_u)^\top$,
$P=A\,B$ is a quotient of polynomials,
\begin{equation}
P=A\,B=\frac{-2\,f_v\,\big(ug-vf\big)\,K}{(\lambda_2-\lambda_1)^2},
\qquad
K=-S_1^{\,2}-M_3\,S_2,\qquad M_3:=e^2-ae+bd-dk .
\label{eq:P}
\end{equation}
The computation behind \eqref{eq:P} is long and tedious; it was carried out by
hand and is omitted here.
The denominator in \eqref{eq:P} is the discriminant of $J$,
$(\lambda_2-\lambda_1)^2=(f_u-g_v)^2+4f_vg_u$, a binary quadratic form that is
positive where the system is strictly hyperbolic, so $P$ vanishes exactly where
the numerator does.

The middle factor of the numerator has a meaning of its own. For a homogeneous quadratic flux
Euler's identity gives $J(U)U=2F(U)$, so $ug-vf=\det[U,\,F(U)]$ vanishes
exactly on the rays along which $U$ is itself an eigenvector of $J(U)$. Each
such ray is an integral curve, and together they form the Hugoniot locus of
the origin, $\{U:\ F(U)=sU\}$, which in the normal form of
Sec.~\ref{sec:normalform} is the cubic of \cite[Appendix]{IMPT88}. The same cubic is the locus where the eigendirection field of the homogeneous model is singular (sorting those singularities into nodes and saddles is Palmeira's classification, recalled in \cite[Sect.~3]{IMPP92}), and its roots are the asymptotic angles of \cite[Eq.~(3.15)]{IMP90}. A direct computation confirms that the three descriptions give one cubic when the flux is homogeneous; once linear terms are added they separate. The factor
vanishes identically only in the degenerate case $F(U)\parallel U$ for every
$U$, in which $P\equiv0$.

That one of $A$ and $B$ vanishes on each such ray is geometric. For a quadratic
flux $F(U+X)=F(U)+J(U)X+F(X)$, so that $J(U)X=J(X)U$ for all $U$ and $X$, and
$\pa_XJ=J(X)$ because $J$ is linear. Let $U=\bbr_k$ on a ray, and let $\bbl_k$ be the
left eigenvector. For $j\neq k$,
\[
\nabla\lambda_k\cdot\bbr_j
=\frac{\bbl_k\,J(\bbr_j)\,\bbr_k}{\bbl_k\cdot\bbr_k}
=\frac{\bbl_k\,J(U)\,\bbr_j}{\bbl_k\cdot\bbr_k}
=\lambda_j\,\frac{\bbl_k\cdot\bbr_j}{\bbl_k\cdot\bbr_k}=0 ,
\]
since $\bbl_k\cdot\bbr_j=0$ for $j\neq k$. Thus the speed of the family whose
eigenvector is the ray is stationary along the curves of the other family.
The factor $f_v$, on the other hand, is an artifact of the normalization, since
on $f_v=0$ one of the $\bbr_i$ vanishes.
\end{remark}

Theorem~\ref{thm:classif} says which systems qualify; let us see what
qualifying costs.

\begin{corollary}
\label{cor:nogo}
In case $(v)$ one has $\lambda_1=\alpha w$ and $\lambda_2=\beta z$,
which coincide on the line $\alpha w=\beta z$. Hence every quadratic system
admitting a closed-form implicit solution fails to be strictly hyperbolic on a
full line through the origin. In particular, no homogeneous quadratic flux
that is strictly hyperbolic away from isolated umbilic points (and so, in the
Schaeffer--Shearer normal form \cite{SS87}, no model off the curve $a=1+b^2$ of
Corollary~\ref{cor:normalform}) admits such a representation.
\end{corollary}

The conclusion is worth restating in the terms a practitioner would use.
Asking for a closed-form implicit solution is asking that each phase be
transportable at a speed depending on its own saturation alone. By
Corollary~\ref{cor:nogo} the only quadratic models for which this is possible
are those whose two characteristic speeds agree along a whole line of states,
so that the two phases cannot be told apart there; a model with an isolated
umbilic point, which is the generic situation and the one the
Schaeffer--Shearer classification was built for, never qualifies. Note that
the obstruction is global, since a decoupled Burgers pair is strictly
hyperbolic on any domain $D$ avoiding the line $\alpha w=\beta z$, so that
Corollary~\ref{cor:nogo} excludes the models with an isolated umbilic point,
among them the Schaeffer--Shearer normal form off the curve $a=1+b^2$, rather
than asserting a pointwise impossibility.

\begin{remark}[the quadratic Corey model]
\label{rem:corey}
Corey fluxes (Remark~\ref{rem:oil}) are rational functions, not quadratic
polynomials, so Theorem~\ref{thm:classif} does not apply, but the criterion
$A\equiv B\equiv0$ does, and since $A$ is continuous one state with $A\neq0$
suffices to show that $A\not\equiv0$ holds. Let
$g=(v^2/\mug)\big/\big(u^2/\muw+v^2/\mug+(1-u-v)^2/\muo\big)$ be the gas
fractional flow, put $\muwg:=\muw+\mug$, $\muwgo:=\muwg+\muo$ and
$\mathbf v:=(\muw,\,\mug)$, and parametrize the line $BO$ of
Remark~\ref{rem:oil} by $U(s)=s\,\mathbf v$, with $O$ at $s=0$ and the umbilic
point at $s=1/\muwgo$. With
$D(s)=s^2\muwg+(1-s\muwg)^2/\muo$ one has
$F(U(s))=\mathfrak f(s)\,\mathbf v$ and $\mathfrak f=s^2/D$, so that the line is
invariant under the flux, and differentiating with respect to $s$ gives
$J\,\mathbf v=\mathfrak f'(s)\,\mathbf v$. Hence the line is an
integral curve of one of the eigenvector fields, and $\mathbf v$ is an
eigenvector. Moreover one finds $f_u+g_v=(4sD-s^2D_s)/D^2$ on
the line, with $D_s=\dd D/\dd s$, and $\mathfrak f'=(2sD-s^2D_s)/D^2$, so the other eigenvalue is
$\mathfrak e(s)=2s/D$, and their difference is
\[
\mathfrak f'-\mathfrak e=\frac{2s^2\muwg\,(1-\muwgo s)}{\muo D^2}>0
\qquad\text{for }0<s<1/\muwgo
\]
(the inequality reverses for $1/\muwgo<s<1/\muwg$, that is, between the
umbilic point and $B$).
Thus on the open segment $OU$ the vector $\mathbf v$ is the eigenvector of the fast family, $\mathfrak e=\lambda_1$
and $\mathfrak f'=\lambda_2$, and the system is strictly hyperbolic on the open
segment. A computation gives
$\mathfrak e'(s)=2(1-\muwg\muwgo s^2)/(\muo D^2)$, positive for
$0\le s\le1/\muwgo$ because $\muwg\muwgo s^2$ is at most $\muwg/\muwgo$, which is less than $1$. Since
$f_v\neq0$ on the open segment (by Remark~\ref{rem:oil} it vanishes only for
$u=0$ or $v\muo=(1-u-v)\mug$, and on $OU$ the second condition reads
$s\muo=1-s\muwg$, that is, $s=1/\muwgo$, the umbilic point), $\bbr_2$ is a
nonzero multiple of $\mathbf v$, and $A=\nabla\lambda_1\cdot\bbr_2$ is a nonzero
multiple of $\mathfrak e'$. Hence $A\neq0$ at every point of the open
segment $OU$, for all $\muw,\,\mug,\,\muo>0$, so \eqref{eq:decoupling} fails
for the Corey model and the substitution $X_i=x-t\lambda_i$ is not legitimate
for arbitrary data, for which the exact representation is \eqref{eq:exactrep},
with feet defined by an ODE. For data supported on the line $BO$, however, which
is an integral curve on each side of the umbilic point (where $f_v\neq0$), the
representation of Theorem~\ref{thm:exact} applies to any Cauchy datum restricted
to it.
\end{remark}

\begin{remark}[what is not in the literature]
\label{rem:literature}
The condition $A\equiv B\equiv0$ appears to have no standard name. Its
neighbors are close but distinct. The Temple class \cite{Tem83} is defined by
coinciding shock and rarefaction curves, equivalently by rectilinear
rarefaction curves or by linear degeneracy of the family. The property can
also hold along a single curve instead of on the whole phase space. Following
the invariant submanifolds of Temple \cite{Tem83b}, a \emph{contact manifold}
is a curve on which the rarefaction curve of one family coincides with part of
the Hugoniot locus and the characteristic speed is constant, and in a model of
bidisperse suspensions such curves are decisive for the structure of the
Riemann solutions \cite[Def.~4 and Lemma~2]{BC23}. The
semi-Hamiltonian condition of Tsarev involves three distinct indices and is
therefore vacuous for $n=2$. Ben-Artzi and Li \cite[Def.~21]{BAL07} call a system
weakly coupled when only a partial set of Riemann invariants exists, the
missing ones being coupled in a way that still allows a diagonalized
treatment; this concerns systems of $m\ge3$ equations and is not our
condition. Neither \cite{IMPT88} nor
\cite{BGR05} contains the computation of this section, since the former
characterizes the Region II/III boundary of Corollary~\ref{cor:normalform}
through the lines $\lambda=0$ and the Hugoniot locus of the origin, never
through decoupling or implicit solutions, and the latter treats the projective
geometry of the rarefaction foliation without reference to Riemann invariants
or to the Cauchy problem.
\end{remark}

\subsection{The class in the Schaeffer--Shearer normal form}
\label{sec:normalform}

The criterion becomes a single equation in the normal form of
\cite{IMPT88,SS87}, namely
\begin{equation}
u_t+\tfrac12\big(au^2+2buv+v^2\big)_x=0,\qquad
v_t+\tfrac12\big(bu^2+2uv\big)_x=0,
\label{eq:SSnormal}
\end{equation}
to which every hyperbolic quadratic system is reduced by a linear change of
dependent variables. In the notation of \eqref{eq:quadflux} this is
$(a,\,b,\,c,\,d,\,e,\,k)=(a,\,b,\,1,\,b,\,1,\,0)$, so that $f_v=g_u=bu+v$
identically. Two simplifications follow, and both are peculiar to the normal
form.

First, $S_1=f_v\wdg g_u\equiv0$ for every $(a,\,b)$, since the normal form is
automatically a gradient system, $J=\operatorname{Hess}\Phi$, and the
classification is governed by the single scalar
\[
S_2=f_v\wdg(g_v-f_u)=a-b^2-1,\qquad\text{with}\qquad M_3=-S_2 .
\]
Second, the rows of $\mathcal M$ become $(b,\,1)$, $(a-1,\,b)$ and
$(-b,\,-1)$, so the third is minus the first and rank $\mathcal M\le1$ reduces
to the vanishing of one determinant,
\[
\det\begin{pmatrix} b & 1\\ a-1 & b\end{pmatrix}=b^2-a+1 .
\]
The step in the proof of Theorem~\ref{thm:classif} where the third minor is
deduced from the other two is not needed here. Since $c=1$, one has $f_v\not\equiv0$ for every $(a,\,b)$, so the
triangular branch of that proof cannot occur either.

\begin{corollary}
\label{cor:normalform}
The system \eqref{eq:SSnormal} satisfies $A\equiv B\equiv0$ if and only if
\[
a=1+b^2
\]
and this is precisely the boundary between Regions II and III of the
Schaeffer--Shearer classification \cite{SS87}, as stated in
\cite[p.~1011]{IMPT88}. On that curve the eigendirections are the constant
directions of slope
\[
m_\pm=\frac{-b\pm\sqrt{b^2+4}}{2},
\]
the roots of $m^2+bm-1=0$, and the system is a pair of decoupled Burgers
equations in the corresponding linear coordinates.
\end{corollary}

\begin{figure}[htbp]
\centering
\begin{picture}(220, 185)
\put(0, 0){\includegraphics[width=0.55\linewidth]{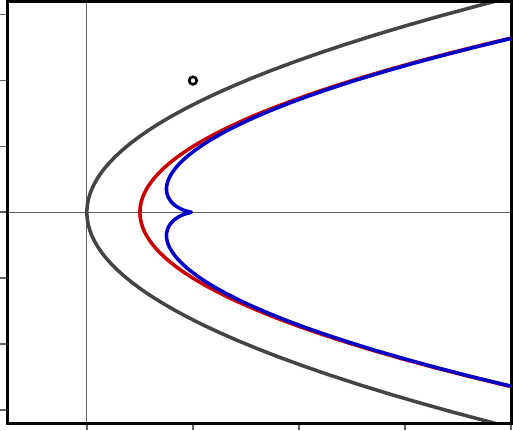}}

\small
\put(209.2, 83.1){$a$}
\put(41.1, 169.8){$b$}

\footnotesize
\put(37.2, -0.9){\makebox(0,0)[t]{$0$}}
\put(82.5, -0.9){\makebox(0,0)[t]{$2$}}
\put(128.3, -0.9){\makebox(0,0)[t]{$4$}}
\put(173.6, -0.9){\makebox(0,0)[t]{$6$}}
\put(219.3, -0.9){\makebox(0,0)[t]{$8$}}
\put(-0.9, 178.6){\makebox(0,0)[r]{$3$}}
\put(-0.9, 150.4){\makebox(0,0)[r]{$2$}}
\put(-0.9, 122.1){\makebox(0,0)[r]{$1$}}
\put(-0.9, 93.8){\makebox(0,0)[r]{$0$}}
\put(-0.9, 65.4){\makebox(0,0)[r]{$-1$}}
\put(-0.9, 37.2){\makebox(0,0)[r]{$-2$}}
\put(-0.9, 8.8){\makebox(0,0)[r]{$-3$}}
\end{picture}
\caption{The $(a,\,b)$ plane. The gray curve $a=\frac34b^2$ is the boundary of
Region~I, the red curve $a=1+b^2$ is the boundary between Regions~II and~III, and
the blue curve is the one between Regions~III and~IV of
\cite[p.~1011]{IMPT88}. The dot is the point
$(a,\,b)=(2,\,2)$ of Sec.~\ref{sec:example}.}
\label{fig:ab}
\end{figure}

Fig.~\ref{fig:ab} places the curve $a=1+b^2$ among the boundaries of the
classification of \cite{SS87}, together with the model of
Sec.~\ref{sec:example}.

We have reduced the three descriptions of the curve $a=1+b^2$ (decoupling,
spectral and combinatorial) to the calculation above, which has two immediate
corollaries, worth spelling out because each is recognizable from a
different direction. The first is spectral. For \eqref{eq:SSnormal} the discriminant is
\[
(\operatorname{tr}J)^2-4|J|=\big[(a-1)^2+4b^2\big]u^2+2b(a+3)uv+(b^2+4)v^2,
\]
a binary quadratic form whose own discriminant equals $-16(a-b^2-1)^2$ times a
positive factor, so it is a perfect square exactly when $a=1+b^2$, and then
\[
(\operatorname{tr}J)^2-4|J|=(b^2+4)\,(bu+v)^2=(b^2+4)\,f_v^2 .
\]
This is not independent evidence, since $\lambda_2-\lambda_1=(m_2-m_1)f_v$,
constant eigendirections force the discriminant to be a perfect square at
once, and $(m_2-m_1)^2=b^2+4$ from $m^2+bm-1=0$. What the computation does add
is the geometry. On the curve $a=1+b^2$ the characteristic speeds coincide along the whole
line $f_v=0$, rather than at an isolated point, so that it is a \textbf{coincidence line},
in the sense of \cite[Def.~3.4]{IMPP92}, where the coincidence locus is the
set of states at which an eigenvalue is multiple (see also \cite{RM13}). Since the normal form is symmetric, $J$ is a multiple of
the identity on that line, so its points are umbilic in the sense of
the same definition, the eigenvalue having geometric multiplicity two. The
coincidence locus is the zero set of the discriminant, a quadratic form in
$(u,\,v)$, hence a conic, which for the normal form is degenerate, being generically
the single point $U=0$, the isolated umbilic point, and on the curve
$a=1+b^2$ the double line $f_v^2=0$. These
systems therefore sit exactly at the edge of the framework of
\cite{IMPT88}, which presumes an isolated umbilic point
\cite[p.~1010]{IMPT88}, and that is presumably why our class was never
examined from the present point of view. The coincidence line is, moreover,
the turning-point locus of \eqref{eq:gammaPrime} discussed in
Sec.~\ref{sec:lifting}.

The second is combinatorial, and it is the one piece of external corroboration. In \cite[p.~1011]{IMPT88} the same boundary is characterized by
the requirement that the two lines $\lambda=0$ coincide with two lines of the
Hugoniot locus of the origin other than the $u$-axis. Substituting
$a=1+b^2$, the lines $\lambda=0$, given by $(a-b^2)u^2-buv-v^2=0$, become
$v=m_\pm u$; and $v=m_\pm u$ annihilates the cubic
$-bu^3+(a-2)u^2v+2buv^2+v^3$ of that Hugoniot locus. The curve was thus
identified in \cite{IMPT88} by counting lines, with no reference to implicit
solutions.

\begin{remark}[the quadratic hypothesis is essential]
\label{rem:gamma3}
Isentropic gas dynamics with exponent $3$, in the form of Remark~\ref{rem:gas},
has Riemann invariants $w=\hat u-c$ and $z=\hat u+c$ and eigenvalues
$\lambda_1=w$, $\lambda_2=z$, so it satisfies $A\equiv B\equiv0$ and decouples
into two Burgers equations \cite{CF48}. Yet its Riemann invariants are not
linear in the conserved variables, its eigendirections are not constant, and
it is known not to belong to the Temple class \cite{Gol22}. This is what Lemma~\ref{lem:invar} leads one to expect, since condition $(i)$ is invariant under every diffeomorphism of the dependent variables, whereas $(iv)$ and $(v)$ are invariant only under linear ones, so that the implications $(i)\,\Rightarrow\,(iv)\,\Rightarrow\,(v)$ need a linear structure, which in Theorem~\ref{thm:classif} is the homogeneity of the flux. The theorem is therefore not a special case of a general statement.
\end{remark}

\begin{remark}[relation to the Temple class]
\label{rem:temple}
Since the eigendirections are constant in case $(iv)$, the rarefaction curves
are straight and parallel within each family, hence they coincide with the
shock curves and the class of Theorem~\ref{thm:classif} is contained in the
Temple class \cite{Tem83}. The inclusion is strict in two ways, since Temple allows
non-parallel straight rarefaction curves, and it allows linearly degenerate
families, whereas the decoupled Burgers pair is genuinely nonlinear.
Conversely Remark~\ref{rem:gamma3} shows that neither class contains the other
once the quadratic hypothesis is dropped.
\end{remark}

\begin{remark}[the wave manifold]
\label{rem:wavemanifold}
The classification has a transparent reading in the wave-manifold formalism
of \cite{AEMP10,IMPP92}. There a shock is a pair $(U_-,\,U_+)$ satisfying the Rankine--Hugoniot
condition, and the set of such pairs, with the diagonal blown up by recording
the projective slope $Z$ of the chord $U_+-U_-$, is the wave manifold; the
\emph{characteristic surface}
$\Sigma$ is the set of shocks of zero amplitude, on which $Z$ becomes the
slope of an eigenvector, so that $Z$ restricted to $\Sigma$ is our $\gamma'$. For
quadratic fluxes the equation of the wave manifold does not involve the length
of the chord, which is the midpoint rule of \cite{IMPT88}, and it is
\eqref{eq:gammaPrime} evaluated at the midpoint of the chord. In the notation
of \cite[(A.8)--(A.9)]{AEMP10} its homogeneous part reads
\[
G=\bar V\,\big(1-b_2Z-Z^2\big)-b_1Z\,\bar U=0 ,
\]
so the eigenvector slopes are independent of the state exactly when $b_1=0$,
which is therefore the class of Theorem~\ref{thm:classif} in that normal form,
as one also checks from $S_1=0$, $S_2=b_1$ and $M_3=-b_1$. On that class
$G=-\bar V\,(Z^2+b_2Z-1)$ factors, its gradient vanishes on the set
$\{\bar V=0,\ Z=Z_\pm\}$, and the wave manifold fails to be a manifold there.
The regularity assumed in the proof of \cite[Thm.~5.1]{AEMP10}, that $\nabla G$
never vanishes for quadratic fluxes, holds for generic lower-order terms but
fails precisely on our class when only the homogeneous part is kept, which
is the local model at an umbilic point. In the language of
\cite[Sect.~3]{IMPP92}, where $\Sigma$ is a ruled surface whose rules are
the fixed-eigenvector curves, straight lines in state space, our class is
the one in which two of these rules fill the whole plane.
\end{remark}

%==============================================================================
\section{The quadratic model and its shocks}
\label{sec:example}

Let us now run the whole argument on one model where every quantity can be
computed. The model is the normal form with $(a,\,b)=(2,\,2)$, which
Sec.~\ref{sec:normalform} has already excluded from the decoupled class, and
the purpose of this section is to see what that exclusion looks like in
practice. In Sec.~\ref{sec:model} we place the model in the Schaeffer--Shearer
classification; in Sec.~\ref{sec:transitional} we recall its crossing discontinuity; in
Sec.~\ref{sec:gnl} we read that discontinuity through the invariant coefficients
$\nabla\lambda_i\cdot\bbr_i$, where the component speeds of
Sec.~\ref{sec:implicit} could not be trusted; in Sec.~\ref{sec:experiments}
we turn to Cauchy data on the diagonal; and in Sec.~\ref{sec:exactshock} we
use the structure found there to write down an exact solution with a
stationary transitional shock.

\subsection{The model}
\label{sec:model}

Following \cite{SS87}, studied also in \cite{EetA25,WePl07}, take
\begin{equation}
f(U)=u^2+2uv+\tfrac12v^2,\qquad g(U)=u^2+uv,\qquad
J(U)=\begin{bmatrix}2u+2v & 2u+v\\ 2u+v & u\end{bmatrix}.
\label{eq:quadmodel}
\end{equation}
Three structural facts organize everything that follows.

\begin{enumerate}
\item \textbf{It is the normal form \eqref{eq:SSnormal} with
$(a,\,b)=(2,\,2)$.} Indeed $\frac12(2u^2+2\cdot2\,uv+v^2)
=u^2+2uv+\frac12v^2$ and $\frac12(2u^2+2uv)=u^2+uv$. In particular it is a
gradient system, $J=\operatorname{Hess}\Phi$ with
$\Phi=\frac16\big(2u^3+6u^2v+3uv^2\big)$ (automatic for the normal form,
since there $f_v=g_u$), so the eigenvectors are orthogonal,
$\det[\bbr_1,\,\bbr_2]=\pm1$ with unit normalization, and \eqref{eq:JacAB}
reduces to $\Jac=\mp(\lambda_2-\lambda_1)\kappa_1\kappa_2$.
\item \textbf{Strict hyperbolicity off the origin, and Region I.} The
discriminant is $(\operatorname{tr}J)^2-4|J|=17u^2+20uv+8v^2$, whose own
discriminant is $20^2-4\cdot17\cdot8=-144<0$ with positive leading
coefficient, which makes it positive definite, so $\lambda_1<\lambda_2$ for all
$U\neq0$, with a single umbilic point at the origin. Moreover
$4a-3b^2=-4<0$, so $(2,\,2)$ lies in Region I of \cite{SS87} (see Fig.~\ref{fig:ab}), where
$\lambda=0$ only at the origin; and since $b=2\neq0$ it is not one of the
symmetric systems $b=0$ for which \cite{IMPT88} and its sequels solve the
Riemann problem.
\item \textbf{No closed form.} By Corollary~\ref{cor:normalform} the criterion
is $a=1+b^2$, whereas here the inequality $a=2\neq5=1+b^2$ holds;
equivalently, $S_2=a-b^2-1=-3$ does not vanish. The model is excluded, consistently with item~2.
\end{enumerate}

Furthermore, \eqref{eq:gammaPrime} is homogeneous of degree $0$ for this flux,
so the substitution $v=u\zeta$ separates variables and the integral curves
are given by the quadrature
\begin{equation}
\frac{\dd u}{u}=\frac{2(2+\zeta)\,\dd\zeta}
{-2\zeta^2-6\zeta-1\pm\sqrt{8\zeta^2+20\zeta+17}} ,
\label{eq:quadrature}
\end{equation}
and the radicand being quadratic, the integral is elementary.
Theorem~\ref{thm:exact} is therefore fully constructive for this model.

With $\Delta:=17u^2+20uv+8v^2$ one has $f_u+g_v=3u+2v$,
$\lambda_{1,2}=\frac12\big(3u+2v\mp\sqrt\Delta\big)$ and
$\bbr_{1,2}=\big(2u+v,\ \tfrac12(-u-2v\mp\sqrt\Delta)\big)^\top$, and the
integral curves of the $i$-th family are those of $\dot U=\bbr_i(U)$; they are
drawn in the background of Fig.~\ref{fig:hugoniot}.

\subsection{A crossing discontinuity}
\label{sec:transitional}

\begin{figure}[htbp]
\centering
\small
\begin{picture}(300, 140)
\put(0, 0){\includegraphics[width=0.75\linewidth]{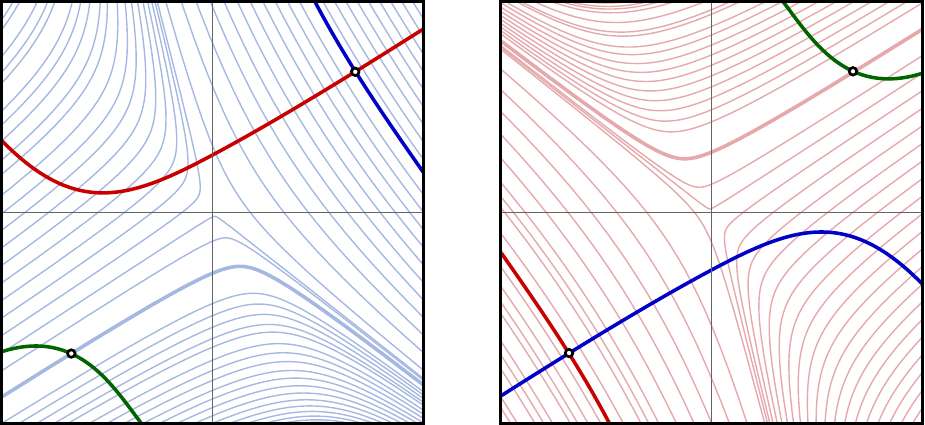}}

\put(289,  62){$u$}
\put(127,  62){$u$}
\put( 66, 130){$v$}
\put(222, 130){$v$}

\put(3, 128){$\mathcal{H}(U_+)$}
\put(111, 111){\makebox(0,0)[rb]{$U_+$}}
\put( 25,  25){$U_-$}

\put(165, 128){$\mathcal{H}(U_-)$}
\put(261, 109){$U_+$}
\put(189,  17){$U_-$}
\end{picture}
\caption{Hugoniot loci for $U_+=(1,\,1)$ on the left and $U_-=(-1,\,-1)$ on
the right, in the domain $u,\,v\in[-3/2,\,3/2]$. Thin straight lines denote
$u=0$ or $v=0$. Colors indicate the shock family, with blue for the slow family,
red for the fast one and green for the crossing discontinuity. Faint curves in the
background are integral curves of the slow family (left) and of the fast family
(right).}
\label{fig:hugoniot}
\end{figure}

Take $U_+=(1,\,1)$ and $U_-=(-1,\,-1)$. Recall that the Hugoniot locus $\mathcal H(U_\pm)$ (met for the origin in
Remark~\ref{rem:P}) is the set of states $U$ for which
the Rankine--Hugoniot condition
\begin{equation}
f(U)-f(U_\pm)=\sigma\,(u-u_\pm),\qquad
g(U)-g(U_\pm)=\sigma\,(v-v_\pm)
\label{eq:RH}
\end{equation}
holds for some $\sigma$; eliminating $\sigma$,
\[
\mathcal H(U_\pm)=\Big\{U\in\bbR^2 \;\Big|\;
\Big(u^2+2uv+\tfrac12v^2-\tfrac72\Big)(v\mp1)=(u^2+uv-2)(u\mp1)\Big\}.
\]
Here $U_\pm$ are double zeros of $\mathcal H(U_\pm)$ and simple zeros of
$\mathcal H(U_\mp)$; see Fig.~\ref{fig:hugoniot}. Notice that $U_-$ lies on
the non-local branch of $\mathcal H(U_+)$ and conversely, whereas in
petroleum models such as those of Remark~\ref{rem:oil} it is common that if $U$ is non-local
in $\mathcal H(V)$ then $V$ lies on a local branch of $\mathcal H(U)$. The
shock speed is symmetric in $U_+$ and $U_-$, and $\sigma=0$ holds, while
\[
\lambda_{1,2}(U_+)=\frac{5\pm\sqrt{45}}{2},\qquad
\lambda_{1,2}(U_-)=\frac{-5\pm\sqrt{45}}{2},
\]
so $\lambda_1(U_\pm)<0<\lambda_2(U_\pm)$ holds. The discontinuity therefore
crosses the characteristic speeds, $\lambda_1<\sigma<\lambda_2$ at both of its
states, and it is a Lax shock of neither family, and in the terminology of
\cite{IMP90} it is a \textbf{crossing discontinuity}. When a crossing
discontinuity has a viscous profile (a saddle--saddle connection) it is called
a \textbf{transitional shock wave} \cite{IMP90,LetA24}, a term that
\cite{LetA24} takes as a synonym of undercompressive shock (see also
\cite{Loz18} for three-phase flow). Whether the jump between $\pm(1,\,1)$ is
transitional is decided in Sec.~\ref{sec:exactshock}. The stability and
$L^p$-decay of shock waves is delicate even in the scalar case \cite{Pla11}.

\subsection{Near linear degeneracy at the crossing discontinuity}
\label{sec:gnl}

By Remark~\ref{rem:noninvariance} the component speeds $\lambda$ and $\mu$
cannot carry an admissibility statement. The invariant quantities to look at
are the genuine-nonlinearity coefficients $\nabla\lambda_i\cdot\bbr_i$, and
for this model they are available in closed form. Recall that
$\Delta=17u^2+20uv+8v^2$ is the discriminant. With the unnormalized
eigenvectors $\bbr_i=(f_v,\,\lambda_i-f_u)^\top$ of the proof of
Theorem~\ref{thm:classif},
\begin{equation}
\nabla\lambda_1\cdot\bbr_1
=-\,\frac{(2u+v)\big(23u+10v-5\sqrt{\Delta}\big)}{2\sqrt{\Delta}},
\qquad
\nabla\lambda_2\cdot\bbr_2
=\frac{(2u+v)\big(23u+10v+5\sqrt{\Delta}\big)}{2\sqrt{\Delta}} .
\label{eq:gnl}
\end{equation}
Both carry the factor $2u+v=f_v$, the turning-point locus of
\eqref{eq:gammaPrime}, and it plays a different role in each family. On
$f_v=0$ the matrix $J$ is diagonal, since the normal form has $f_v=g_u$; one
eigenvector is $(0,\,1)^\top$ and the other is horizontal. For the family with
the horizontal eigenvector, the unnormalized $\bbr_i=(f_v,\,\lambda_i-f_u)^\top$
itself vanishes on $f_v=0$, and the factor cancels on normalizing. For the
family with the vertical eigenvector the zero is genuine, because
$\nabla\lambda\cdot(0,\,1)^\top=\pa_v g_v=g_{vv}$, and $g_{vv}=k$ vanishes
identically in the normal form. So the line $f_v=0$ is always an
inflection line of the normal form, for whichever family has its integral
curve vertical there; in the present model that is the fast family for $u>0$
and the slow one for $u<0$. (The line $f_v=0$, that is $v=-bu$, lies in the
Hugoniot locus of the origin exactly when $b(1-a+b^2)=0$, as substitution
into the cubic of Sec.~\ref{sec:normalform} shows; the case $b=0$ is consistent
with the observation of \cite[p.~1011]{IMPT88} that for symmetric systems the
lines where genuine nonlinearity fails lie in that locus, and the case $a=1+b^2$ is
the curve of Corollary~\ref{cor:normalform}.) Discarding the factor,
\[
\nabla\lambda_1\cdot\bbr_1=0
\iff 23u+10v=5\sqrt{\Delta}
\iff 25v^2+10uv-26u^2=0
\quad\text{with}\quad 23u+10v>0 ,
\]
and the quadratic factors into two straight lines through the origin,
$v=\zeta_\pm\,u$ with
\begin{equation}
\zeta_\pm=\frac{-1\pm3\sqrt3}{5},
\label{eq:inflection}
\end{equation}
the two roots of $25\zeta^2+10\zeta-26=0$; the subscript records which root is
meant, $\zeta_+>0>\zeta_-$. For $u>0$ these are inflection lines of the slow
family, and the fast family degenerates only on $f_v=0$, since
$23u+10v+5\sqrt\Delta>0$ there. The picture for $u<0$ follows from the
symmetry $U\mapsto-U$ (compare the two panels of Fig.~\ref{fig:hugoniot}). Because
$J$ is linear, $J(-U)=-J(U)$, so
$\lambda_1(-U)=-\lambda_2(U)$ and $\bbr_1(-U)=\bbr_2(U)$, which gives
$\nabla\lambda_1\cdot\bbr_1(-U)=\nabla\lambda_2\cdot\bbr_2(U)$. The map
$U\mapsto-U$ therefore exchanges the two families. The inflection locus of the
model is thus three straight lines through the origin, $v=-2u$ and
$v=\zeta_\pm u$, each of them an inflection line of one family on one side of
the origin and of the other family on the other side.

Let us now place the states of the discontinuity on this picture. The states
$U_\pm=\pm(1,\,1)$ of Sec.~\ref{sec:transitional} lie on the line $v=u$, which
makes an angle of about $5.0^\circ$ with the inflection line $v=\zeta_+u$ (slope
$1$ against $\zeta_+\approx0.8392$), and they lie on opposite sides of the
origin. By the symmetry just described, whatever
happens at one of them happens at the other with the families exchanged, so
it suffices to compute at $U_+$. There the slow coefficient $\nabla\lambda_1\cdot\bbr_1$, computed with unit eigenvectors of positive first component in the Euclidean norm of $(u,\,v)$, is $0.0212$ against $4.219$ for the fast coefficient $\nabla\lambda_2\cdot\bbr_2$, a ratio of about one to two hundred (the value of such a coefficient, unlike its vanishing, depends on the normalization, see Lemma~\ref{lem:invar}); at $U_-=-U_+$ the same two numbers appear with the
families exchanged. Hence at each end of the discontinuity one
family is nearly linearly degenerate (the slow family at $U_+$ and the fast
family at $U_-$). Unlike the comparison of $\lambda$ with $\mu$ in
Remark~\ref{rem:noninvariance}, this is a statement about the system and not
about a choice of variables.

\begin{remark}[an observation]
\label{rem:neardeg}
Near linear degeneracy of a family means that its characteristic speed is
almost constant along its rarefaction curves, so near that state the family
behaves almost like a contact field, so that weak waves of that family barely steepen
or spread, since the rate at which they do so is governed by
$\nabla\lambda_i\cdot\bbr_i$. In this example the nearly degenerate family is,
at each end of the discontinuity, also the one whose characteristic speed is
closer to the shock speed $\sigma=0$, since at $U_+$ one has
$\lambda_1=(5-\sqrt{45})/2\approx-0.85$ against $\lambda_2\approx5.85$, and at
$U_-$ the families are exchanged, with $\lambda_2\approx0.85$ against
$\lambda_1\approx-5.85$. So on each side the family that almost travels with
the discontinuity is also the one that almost fails to be genuinely nonlinear.
The same happens, with the values of Sec.~\ref{sec:exactshock}, at the states
of the exact transitional shock below. I have not found this coincidence in
\cite{IMP90}, whose Eq.~(3.15) and Cor.~3.7 are what Secs.~\ref{sec:normalform}
and~\ref{sec:exactshock} use, and I do not know whether it is a feature of
transitional shocks in general or an accident of these examples. I have no
argument for it, and record it only as an observation.
\end{remark}

\subsection{Cauchy data on the diagonal}
\label{sec:experiments}

Consider
\[
U_{0,\pm}(x)=\pm\Big(\tfrac{2}{\pi}\arctan x,\ \tfrac{2}{\pi}\arctan x\Big),
\]
so that $U_{0,+}(x)\to U_\pm$ as $x\to\pm\infty$ and $U_{0,-}(x)\to U_\pm$ as
$x\to\mp\infty$, reminiscent of shocks seen through viscous profiles
\cite{CetA16}. In both cases the datum is the segment joining $(-1,\,-1)$ to
$(1,\,1)$, so $\gamma(u)=u$ holds with $\gamma'\equiv1$, and \eqref{eq:lammu} gives
$\lambda(u)=7u$ and $\mu(v)=4v$.

\begin{figure}[htbp]
\centering
\small
\begin{picture}(300, 138)
\put(0, 0){\includegraphics[width=0.75\linewidth]{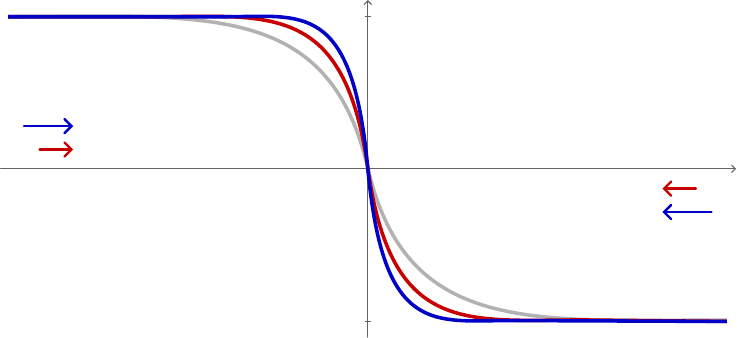}}

\put(290, 60){$x$}
\put(140, 127.5){$1$}
\put(134,   3.5){$-1$}

\put(283,  10){$U_-$}
\put(  4, 120){$U_+$}

{\footnotesize
\put(264, 47){\makebox(0,0)[rb]{$\lambda(u) \;=\; 7u$}}
\put(264, 59){\makebox(0,0)[rb]{$\mu(v) \;=\; 4v$}}}

{\footnotesize
\put(33, 85){$\lambda(u) \;=\; 7u$}
\put(33, 73){$\mu(v) \;=\; 4v$}}
\end{picture}
\caption{First-order predictor \eqref{eq:ansatz} for the datum $U_{0,-}$
joining $U_+$ to $U_-$ (thin gray curve). The $u$ component (blue) advances at
$\lambda(u)$ and the $v$ component (red) at $\mu(v)$. Since
$\lambda\neq\mu$ the two separate, and the separation is the quantity
estimated in Proposition~\ref{prop:error}; this is a prediction of incipient
steepening, not a solution of \eqref{eq:system}.}
\label{fig:atanShock}
\end{figure}

\begin{figure}[htbp]
\centering
\small
\begin{picture}(300, 138)
\put(0, 0){\includegraphics[width=0.75\linewidth]{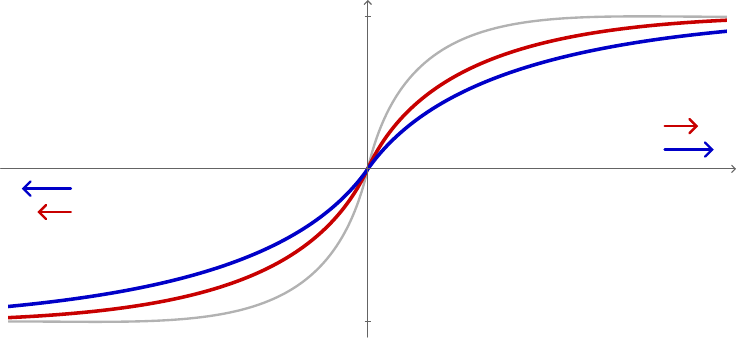}}

\put(290, 60){$x$}
\put(140, 127.5){$1$}
\put(134,   3.5){$-1$}

\put(283, 112){$U_+$}
\put(  4,  17){$U_-$}

{\footnotesize
\put(264, 73){\makebox(0,0)[rb]{$\lambda(u) \;=\; 7u$}}
\put(264, 85){\makebox(0,0)[rb]{$\mu(v) \;=\; 4v$}}}

{\footnotesize
\put(33, 59){$\lambda(u) \;=\; 7u$}
\put(33, 47){$\mu(v) \;=\; 4v$}}
\end{picture}
\caption{The same predictor for the datum $U_{0,+}$ joining $U_-$ to $U_+$.
Both components spread, the $u$ profile faster; no steepening is predicted.
Same caveat as in Fig.~\ref{fig:atanShock}.}
\label{fig:atanRar}
\end{figure}

Notice that $J(U)(1,\,1)^\top=(7u,\,4u)^\top$ is not parallel to
$(1,\,1)^\top$ except at the origin, so the datum does not lie on an integral
curve of an eigenvector field, Theorem~\ref{thm:exact} does not apply, and \eqref{eq:ansatz} is only
the first-order predictor of Proposition~\ref{prop:error}, with $\Gamma_t$
leaving the diagonal at the rate given there. Within that accuracy, the $u$
profile steepens faster than the $v$ profile for $U_{0,-}$
(Fig.~\ref{fig:atanShock}) and spreads faster for $U_{0,+}$
(Fig.~\ref{fig:atanRar}).

Lax--Friedrichs simulations of RP$(U_+,\,U_-)$ and RP$(U_-,\,U_+)$ confirm the
qualitative picture. The former produces a zero-speed shock with left and right states
\[
L\approx(1.1548,\;0.7569)\approx-R,
\]
which lie on the line $v=mu$ with $m\approx0.6554$ (the slope of $L$ is $0.6554$
to four digits). The latter also produces a zero-speed shock, now between
$\pm(0.10,\,-0.40)$ on the line with $m\approx-3.8662$, the third root of
\eqref{eq:cubicm}, flanked by waves of the two families. Remark~\ref{rem:RP}
explains why such states appear.

\begin{figure}[htbp]
\centering
\small
\begin{picture}(380, 85)
\put(0, 0){\includegraphics[width=0.45\linewidth]{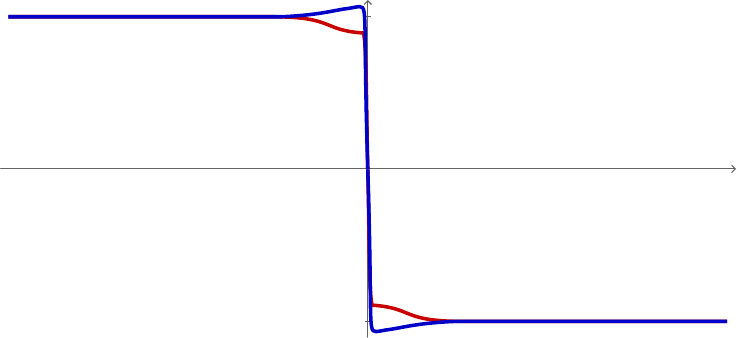}}
\put(195, 0){\includegraphics[width=0.45\linewidth]{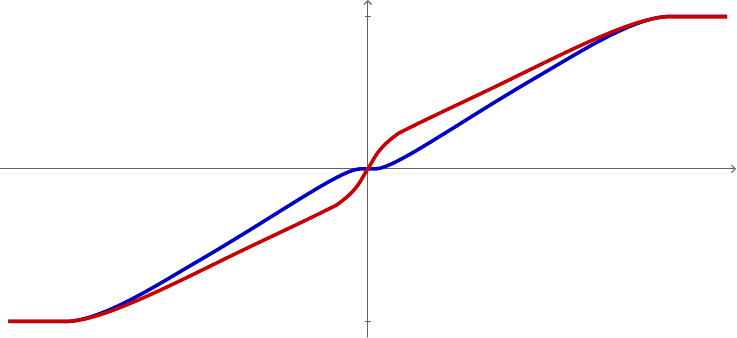}}

\put(170, 33){$x$}
\put(165,  7){$U_-$}
\put(  3, 69){$U_+$}

\put(365, 33){$x$}
\put(360, 69){$U_+$}
\put(198,  7){$U_-$}
\end{picture}
\caption{Lax--Friedrichs simulations for the arctangent profiles (not for
Riemann data), at the final time $t=2$. In the left panel $U_{0,-}$ develops a zero-speed shock flanked by
rarefactions, whose states differ from $U_\pm$, so the crossing discontinuity from $U_+$ to
$U_-$ does not persist, and in the right panel $U_{0,+}$ resolves as a rarefaction.}
\label{fig:RP}
\end{figure}

Lax--Friedrichs simulations for the arctangent profiles, and not for Riemann data,
are shown in Fig.~\ref{fig:RP}. The compressive profile $U_{0,-}$ develops a
zero-speed shock flanked by rarefactions, while $U_{0,+}$ resolves as a
rarefaction. The states on the two sides of the shock are not $U_\pm$. They are
opposite, $R=-L$, as the symmetry $U(x,\,t)\mapsto-U(-x,\,t)$ of the system and
of the datum requires, and they lie on the same line $v=mu$ with $m\approx0.6554$.
The slope of $L$ is $0.6554$ for $t\ge1.5$ (it is $0.6558$ at $t=1$). The first component of $L$ is $1.02$ at $t=1$, $1.12$ at $t=4$ and
$1.14$ at $t=8$, so $L$ moves slowly along the line. The number $m$ is not an
accident, since it is a root of a cubic, as we show in
Sec.~\ref{sec:exactshock}, and the data suggest that $L$ tends to the left state of the Riemann problem
above, which lies on the same line.

\subsection{An exact transitional shock}
\label{sec:exactshock}

The simulations select a line through the origin, and Theorem~\ref{thm:exact}
says when data on a line are exactly solvable, so we use the theory backwards,
as an ansatz. Take $\gamma(u)=mu$ with $m$ constant, so that $\gamma'=m$. For the
model \eqref{eq:quadmodel}, \eqref{eq:lammu} gives
\[
\lambda=f_u+mf_v=(2+4m+m^2)\,u,\qquad
\mu=\frac{g_u}{m}+g_v=\frac{2+2m}{m}\,u ,
\]
and by \eqref{eq:equiv} the equality $\lambda=\mu$ holds if and only if
\begin{equation}
m^3+4m^2-2=0 .
\label{eq:cubicm}
\end{equation}
Asking the line to lie in the Hugoniot locus of the origin, $u\,g=v\,f$ at
$(1,\,m)$, gives $1-2m^2-\frac12m^3=0$, which is the same cubic (\emph{cf.}\
Sec.~\ref{sec:normalform}). Its roots are real, $m\approx0.6554$, $-0.7892$ and
$-3.8662$; the positive one is
\[
m=\frac{8\cos\big(\frac13\arccos(-\frac{37}{64})\big)-4}{3}\approx0.655442,
\]
and we take it. On that line $\lambda$ and $\mu$ both equal $cu$, with $c=2+4m+m^2\approx5.0514$, and
$(1,\,m)$ is an eigenvector of $J$ for the eigenvalue $cu$, the fast one for
$u>0$ and the slow one for $u<0$.

Among the states of this line we take the one closest to $(1,\,1)$ in the
Euclidean norm, which keeps the example next to that of
Sec.~\ref{sec:transitional}. It is $U_+=u_*(1,\,m)$ with
$u_*=(1+m)/(1+m^2)\approx1.1580$, that is $U_+\approx(1.1580,\,0.7590)$, and
we put $U_-=-U_+$. As in Sec.~\ref{sec:transitional}, $\sigma=0$ holds for the
pair because $f$ and $g$ are even, and $\lambda_{1,2}(U_+)\approx-0.8575,\,
5.8494$, so $U_\pm$ is again a crossing discontinuity. The datum
\[
U_0(x)=u_0(x)\,(1,\,m)^\top,\qquad u_0(x)=-\,u_*\,\tfrac{2}{\pi}\arctan x,
\]
joins $U_+$ on the left to $U_-$ on the right, lies on a line that is an integral curve, and
by Theorem~\ref{thm:exact} breaks at
$t_b=\big(c\,u_*\tfrac{2}{\pi}\big)^{-1}=\pi/(2cu_*)\approx0.2685$. What happens
afterwards is explicit.

\begin{proposition}
\label{prop:exactshock}
Let $m$, $c$, $u_*$ and $U_0$ be as above. For $0\le t<t_b$ the function
$U=u\,(1,\,m)^\top$ with $u=u_0(x-ctu)$ is a classical solution of
\eqref{eq:system}. For $t>t_b$ let $\xi_s(t)>0$ be the solution of
$\xi_s/\arctan\xi_s=t/t_b$ and set $u_L(t)=u_*\frac{2}{\pi}\arctan\xi_s(t)$.
Let $u(x,\,t)$ be given by the relations
\[
x=\xi-\frac{t}{t_b}\arctan\xi,\qquad u=-\,u_*\frac{2}{\pi}\arctan\xi ,
\]
with $\xi\ge\xi_s(t)$ for $x>0$ and $\xi\le-\xi_s(t)$ for $x<0$. Then
$U=u\,(1,\,m)^\top$ is a weak solution of \eqref{eq:system} with a single
shock, stationary at $x=0$, which joins $u_L(t)(1,\,m)^\top$ to its opposite and
has a viscous profile.
\end{proposition}

\begin{proof}
On the line, $f(u,\,mu)=ku^2$ and $g(u,\,mu)=(1+m)u^2$ with
$k=1+2m+\frac12m^2$, and \eqref{eq:cubicm} is equivalent to $1+m=mk$. Hence
$F(u(1,\,m)^\top)=ku^2(1,\,m)^\top$, and for $v=mu$ both equations of
\eqref{eq:system} reduce to the scalar law $u_t+(ku^2)_x=0$, that is
$u_t+c\,u\,u_x=0$ since $c=2k$. For $t<t_b$ this is
Theorem~\ref{thm:exact}. For $t>t_b$, $u$ is the entropy solution of this
law for the odd and decreasing datum $u_0$, where the characteristics with feet
$|\xi|<\xi_s(t)$ have been absorbed by a single shock, which by oddness is
stationary at $x=0$ with states $\pm u_L(t)$, and the displayed relations are
the remaining characteristics $x=\xi+ctu_0(\xi)$. Since $f$ and $g$ are even,
$F(U_L)=F(-U_L)$ with $U_L=u_L(1,\,m)^\top$, so the Rankine--Hugoniot
conditions hold with $\sigma=0$. Finally a viscous profile $U(x/\varepsilon)$
with $U(\mp\infty)=\pm U_L$ solves $\varepsilon U'=F(U)-F(U_L)$; the line is
invariant, so with $U=w\,(1,\,m)^\top$ one gets
$\varepsilon w'=k(w^2-u_L^2)$, whose solution
$w=-u_L\tanh(ku_Lx/\varepsilon)$ joins $U_L$ to $-U_L$, which completes the
proof.
\end{proof}

The shock forms as in Burgers' equation, since for $t<t_b$ the characteristics
$x=\xi+ctu_0(\xi)$ do not cross, and at $t=t_b$ the slope of $u$ at $x=0$
becomes infinite. The rupture happens at one point only. The characteristic with foot $\xi$ is
absorbed at time $t_b\,\xi/\arctan\xi$, which is larger than $t_b$ for every
$\xi\neq0$, so the shock is born at $x=0$ and at $t=t_b$ with zero amplitude,
at the umbilic state $U=0$, and the rest of the profile remains smooth. Its
amplitude grows like $\frac{\sqrt3}{ct}\sqrt{1-t_b/t}$ just after $t_b$ and
tends to $u_*$, as $u_L=u_*-2/(\pi ct)+\cdots$. For instance $u_L=0.5127$ at
$t=1.2\,t_b$ and $u_L=0.6634$ at $t=1.4\,t_b$. A Lax--Friedrichs
approximation on a grid of $3\cdot10^{4}$ cells over $[-30,\,30]$ agrees with
the solution of the proposition at $t=0.8,\,1.0,\,1.2$ and $1.4$ times $t_b$,
and away from the shock the two differ by a few $10^{-3}$, as expected of a
first-order scheme, and $v=mu$ holds to rounding error.

\begin{remark}[the Riemann problem]
\label{rem:RP}
The Riemann problem evolves as the transitional shock just computed, both in
theory and under the artificial viscosity of the Lax--Friedrichs scheme.
Consider first Riemann data on the line, $L=a\,(1,\,m)^\top$ and $R=-L$ with $a>0$
and $m\approx0.6554$. As in the proof of Proposition~\ref{prop:exactshock}, the
problem reduces to the scalar one $u_t+c\,u\,u_x=0$ with left state $a$ and right
state $-a$, whose entropy solution is the stationary shock, and the viscous
profile found there shows that the jump between $L$ and $R$ is a transitional
shock. Off the line, for instance for RP$(U_+,\,U_-)$ with $U_\pm=\pm(1,\,1)$,
there is no viscous profile for that jump (by \cite[Cor.~3.7]{IMP90}, as
discussed below), and the transitional shock appears between states $\pm L$ on
the line, where $L$ is the point at which the slow wave curve of $U_+$ meets the
line $v=mu$. Integrating \eqref{eq:quadrature} gives
$L\approx(1.15476,\,0.75688)$, and intersecting the Hugoniot locus of $U_+$ with
the line gives the same point to seven digits in $u$, so the choice between the
two curves does not matter here. The waves joining $U_+$ to $L$ and $-L$ to $U_-$
travel at speeds close to $\lambda_1(U_+)\approx-0.854$ and
$\lambda_2(U_-)\approx0.854$, because the corresponding family is nearly linearly
degenerate there (Sec.~\ref{sec:gnl}).

The Lax--Friedrichs scheme has numerical viscosity $\frac{\Delta x^2}{2\Delta t}$
times the identity matrix, which is the identity viscosity assumed throughout this
section, and it reproduces this picture. On $14\,000$ cells over $[-14,\,14]$ and
at $t=2$ it gives $L\approx(1.1548,\,0.7569)$ for RP$(U_+,\,U_-)$, against
$(1.15476,\,0.75688)$ above, while for the datum $\pm U_+$ of
Sec.~\ref{sec:exactshock}, which lies on the line, it returns the two states
unchanged away from the shock layer.
\end{remark}

At $U_\pm$ the slow coefficient $\nabla\lambda_1\cdot\bbr_1$ and the fast one
$\nabla\lambda_2\cdot\bbr_2$, with unit eigenvectors of positive first
component, are $-0.0282$ and $4.2248$ at $U_+$, a ratio of about $0.0067$, so
at each end of this transitional shock one family is again nearly linearly
degenerate, the slow one at $U_+$ and the fast one at $U_-$. The other two
roots of \eqref{eq:cubicm} give lines with the same property $\lambda\equiv\mu$,
and the computation of the proof applies to them similarly, since only $k$ and $c$ change with $m$, and $1+m=mk$ still holds by \eqref{eq:cubicm}. That $u$ and $v$ have
opposite signs there plays no role, since $v=mu$ merely follows $u$; what
matters is which end of the datum lies on the left. The characteristic speed is
$cu$ with $c=2k$, and the datum is compressive when its left end moves to the
right, which for $m\approx0.6554$ and $m\approx-3.8662$ (where $c>0$) means
$u>0$ on the left, and for $m\approx-0.7892$ (where $c<0$) it means $u<0$ on
the left. The opposite orientation gives a rarefaction. Lax--Friedrichs runs
confirm the three cases, and at the states of the shock for $m\approx-0.7892$
one finds again $\lambda_1<0<\lambda_2$ with $\sigma=0$.

The criterion $\lambda_1<0<\lambda_2$ with $\sigma=0$ holds for every pair
$\pm U$, because $f$ and $g$ are even, so by itself it does not guarantee a
transitional shock. On the three lines $v=mu$ with $m^3+4m^2-2=0$ one has in
addition $\lambda\equiv\mu$, so the line is invariant for the equation of the
viscous profile, which reduces to a scalar equation with an explicit solution,
and the jump between $\pm U$ (with the compressive orientation) has a viscous
profile. Off those lines, for instance for $\pm(1,\,1)$, there is no such
profile, by \cite[Cor.~3.7]{IMP90}, since the flux is a gradient and the viscosity is
the identity, so that a crossing discontinuity has a viscous profile if and only if
its two states lie on the same bifurcation line, that is, on one of the three
lines above, and the second is correctly oriented with respect to the first
(which on those lines is the compressive orientation). The simulation then
replaces the jump by a plateau flanked by two fans. The identity viscosity is
the degenerate case of that analysis, since for other viscosity matrices the
straight-line profiles occur at other angles \cite[Thm.~3.5]{IMP90}.

%==============================================================================
\section{Final remarks}
\label{sec:final}

Let us collect the picture and then the limitations. For Riemann data the
solution is self-similar, hence rank one by Proposition~\ref{prop:rank1},
hence its state path is a curve; the EFF is the lifting of that curve and
satisfies Ole\u{\i}nik's criterion \cite{Cas18}. For general Cauchy data the
Jacobian no longer vanishes, and by \eqref{eq:JacAB} the product of wave
strengths $\kappa_1\kappa_2$ measures how far the state path deforms. The implicit
closed-form representation is then exact on simple waves
(Theorem~\ref{thm:exact}) and first-order accurate off them
(Proposition~\ref{prop:error}), while for arbitrary data it exists only in the
thin class of Theorem~\ref{thm:classif}, which excludes the models of interest
by Corollary~\ref{cor:nogo}. What remains valid in general is the ODE-based
representation \eqref{eq:exactrep}.

For applications the useful content is quantitative. The initial propagation
speeds
\[
\lambda(u(x,\,0))=f_u(U_0)+\frac{v_0'(x)}{u_0'(x)}f_v(U_0),
\qquad
\mu(v(x,\,0))=\frac{u_0'(x)}{v_0'(x)}g_u(U_0)+g_v(U_0)
\]
are explicit from the data, and Proposition~\ref{prop:error} bounds the
departure of the state path from $\Gamma$ by $O(t|\lambda-\mu|)$. This is the
form in which the construction supports Lagrangian--Eulerian implementations
\cite{AetA22,AetA23}, and it applies without change to gravitational
multiphase models \cite{Rod23,RMd25}.

Let us also record a heuristic that comes out of Proposition~\ref{prop:error}
and Sec.~\ref{sec:experiments}; we do not claim it as a theorem. Take Cauchy
data joining the two states of a crossing discontinuity, and compare the
initial speeds $\lambda$ and $\mu$ of the two components. If the data make both
components spread, as for $U_{0,+}$ (Fig.~\ref{fig:atanRar}), the predictor
shows no steepening, and the argument indicates that no shock forms from such
data, at least at first. If instead they steepen, as for $U_{0,-}$
(Fig.~\ref{fig:atanShock}), the predictor says that a shock may form, but it
cannot say which one, since it holds for short times only. In the example the
shock that forms is not the discontinuity between $\pm(1,\,1)$ but one between
states $\pm L$ on a line where $\lambda\equiv\mu$ (Sec.~\ref{sec:exactshock}),
and that extra hypothesis is exactly the one under which
Theorem~\ref{thm:exact} makes the solution explicit. Thus the spreading
direction indicates that no shock forms, whereas the compressive direction
needs the extra hypothesis but, with it, leads to analytical results for the
elusive transitional shocks. Since $\lambda$ and $\mu$ depend on the choice of
variables (Remark~\ref{rem:noninvariance}), the comparison is to be read in the
variables of the model at hand; the invariant information is where
$\lambda_{1,2}(U_\pm)$ lie relative to $\sigma$.

Let us be explicit about what is not achieved. We have no extension of
Theorem~\ref{thm:classif} beyond homogeneous quadratic fluxes, and
Remark~\ref{rem:gamma3} shows that the proof given here cannot be adapted
without a new idea, since gas dynamics with exponent $3$ satisfies the criterion
without constant eigendirections. For the quadratic Corey fluxes of
Remark~\ref{rem:oil} the criterion $A\equiv B\equiv0$ fails for every choice of
viscosities (Remark~\ref{rem:corey}), so no closed-form implicit solution
exists for them, and the exact representation is \eqref{eq:exactrep}. And the
predictor \eqref{eq:ansatz} is a statement about
short times only, and nothing here describes the state path of the true solution
once the two component profiles have separated. A natural extension is to
systems with accumulation, $G(U)_t+F(U)_x=0$, which arise in foam flow in
porous media \cite{DetA25,deeT23}.

\section*{Acknowledgments}

This study began during my sabbatical year at the Institute of Mathematics of
UNAM at Ciudad Universitaria, and I thank the Institute for the opportunity to
spend that year at its facilities, for the collaboration of colleagues such as
Antonio Capella, and for the environment in which this work was carried out. I
also acknowledge the use of Claude (Anthropic) as an assistant for \LaTeX\
editing, symbolic and numerical checks, and the preparation of figures; the
mathematical content, and any errors it may contain, are my responsibility.

\section*{Conflicts of interest}

The author does not work for, advise, own shares in, or receive funds from any
organization that could benefit from this article, and has declared no
affiliations other than his research organizations.

%==============================================================================

\end{document}